\documentclass{amsart}
\usepackage[utf8]{inputenc}
\usepackage{amsfonts}
\usepackage{hyperref}
\usepackage{amsmath}
\usepackage{xcolor}
\usepackage{color}
\usepackage{amsthm}
\usepackage{pdflscape}
\usepackage{mathrsfs}
\usepackage{verbatim}
\usepackage{mathtools}
\usepackage{dsfont}
\usepackage{bbm}
\usepackage{amssymb}
\usepackage{geometry}
\usepackage{amstext}
\usepackage{enumerate}
\usepackage{graphicx}
\usepackage{subfig}

\usepackage{geometry}
\usepackage{scalerel,stackengine}
\stackMath
\newcommand\reallywidehat[1]{%
\savestack{\tmpbox}{\stretchto{%
  \scaleto{%
    \scalerel*[\widthof{\ensuremath{#1}}]{\kern-.6pt\bigwedge\kern-.6pt}%
    {\rule[-\textheight/2]{1ex}{\textheight}}
  }{\textheight}%
}{0.5ex}}%
\stackon[1pt]{#1}{\tmpbox}%
}

\DeclareMathOperator{\vol}{vol}
\DeclareMathOperator{\interior}{int}

\hypersetup{colorlinks,%
citecolor=red,%
linkcolor=blue,%
}

\theoremstyle{plain}
\newtheorem*{thm*}{Theorem}
\newtheorem{thm}{Theorem}[]
\newtheorem{cor}{Corollary}[]

\newtheorem{lem}{Lemma}[]
\newtheorem{conj}{Conjecture}[]

\newtheorem{rem}{Remark}[]
\newtheorem{example}{Example}[]
\newtheorem{problem}{Problem}[]
\newtheorem*{observation*}{Observation}

\theoremstyle{definition}

\newcommand{\tvect}[2]{\tiny{
    \begin{pmatrix} #1 \\ #2  \end{pmatrix} }}
\newcommand{\A}{{\bf A}}

\newcommand{\Z}{\mathbb{Z}}
\newcommand{\R}{\mathbb{R}}
\newcommand{\C}{\mathbb{C}}

\newcommand{\Lat}{\mathcal{L}}
\newcommand{\PP}{\mathcal{P}}
\newcommand{\QQ}{\mathcal{Q}}
\newcommand{\BB}{\mathcal{B}}
\newcommand{\Int}{\operatorname{int}}

\newcommand{\blue}[1]{\textcolor{blue}{#1}}

\newcommand{\sinc}{{\rm sinc}}

\title[Sharp inequalities for discrete and continuous multi-tiling]{Sharp inequalities for discrete and continuous multi-tiling, using the Bombieri-Siegel approach.}

\author{Michel Faleiros Martins}
\address{Michel Faleiros Martins, Instituto de Matem\'atica e Estat\'\i stica, Universidade de S\~ao Paulo\\ Rua do Mat\~ao 1010, 05508-090 S\~ao Paulo/SP, Brazil.}
\email{michelfaleiros@gmail.com}

\author{Sinai Robins}
\address{Sinai Robins, Instituto de Matem\'atica e Estat\'\i stica, Universidade de S\~ao Paulo\\ Rua do Mat\~ao 1010, 05508-090 S\~ao Paulo/SP, Brazil.}
\email{srobins@ime.usp.br}  


\makeatletter
\def\@tocline#1#2#3#4#5#6#7{\relax
  \ifnum #1>\c@tocdepth 
  \else
    \par \addpenalty\@secpenalty\addvspace{#2}%
    \begingroup \hyphenpenalty\@M
    \@ifempty{#4}{%
      \@tempdima\csname r@tocindent\number#1\endcsname\relax
    }{%
      \@tempdima#4\relax
    }%
    \parindent\z@ \leftskip#3\relax \advance\leftskip\@tempdima\relax
    \rightskip\@pnumwidth plus4em \parfillskip-\@pnumwidth
    #5\leavevmode\hskip-\@tempdima
      \ifcase #1
       \or\or \hskip 1em \or \hskip 2em \else \hskip 3em \fi%
      #6\nobreak\relax
    \hfill\hbox to\@pnumwidth{\@tocpagenum{#7}}\par
    \nobreak
    \endgroup
  \fi}
\makeatother
\usepackage{setspace}
\AtBeginDocument{%
  \addtolength\abovedisplayskip{0.2\baselineskip}%
  \addtolength\belowdisplayskip{0.2\baselineskip}%
}
\graphicspath{{Images/}}
\begin{document}
\begin{abstract}
Given a finite subset $F$ of integer points in $\Z^d$, it is of interest to seek conditions on $F$ that allow it to multi-tile $\Z^d$ by translations. 
In addition to the continuous multi-tiling results presented here, we also give analogous discrete applications to arithmetic combinatorics.  Namely we give a discretized version of the Bombieri-Siegel formula, namely a finite sum  of discrete covariograms, taken over any finite set of integer points in $\R^d$.
As a consequence, we arrive at a new equivalent condition for multi-tiling $\Z^d$ by translating $F$ with a fixed integer sublattice.

Similar questions related to convex bodies have already been investigated extensively. 
In order to develop lattice sums of the cross covariogram for any two bounded sets 
$A, B\subset \R^d$, we prove a refined continuous version of the classical Bombieri-Siegel formula from the geometry of numbers. To achieve this goal, we use 
a variant of the Poisson Summation formula, adapted for  continuous functions of compact support.

As an application of this refined Bombieri-Siegel formula,
a new characterization of multi-tilings of Euclidean space by translations of a compact set by using a lattice is given.
A further consequence is a spectral formula for the volume of any bounded measurable set.

\end{abstract}
\maketitle
\tableofcontents

\section{Introduction}
We extend some results of Siegel and of Bombieri from the Geometry of numbers by studying the {\bf cross covariogram} of any two compact sets 
$A, B \subset \R^d$, defined by
\begin{equation}\label{def:cross covariogram}
   g_{A, B}(x):= \vol(A\cap (B+x)),
\end{equation}
defined for all $x\in \R^d$. An important special case of 
\eqref{def:cross covariogram}
that comes up here as well is
the {\bf  covariogram} of $A$, denoted by
$ g_{A}(x):= \vol(A\cap (A+x))$.
We define a {\bf body} to be a compact subset of $\R^d$, following the standard conventions of convex geometry.   Given any full-rank lattice $\Lat \subset \R^d$, we study the following sum of the covariagram $g_A$ over the lattice $\Lat$:

\begin{equation} 
\sum_{n \in \Lat} g_A(n) :=
\sum_{n \in \Lat}   \vol(A\cap (A+n)),
\end{equation}
which is always a finite sum due to the compactness of $A$.  The covariogram $g_A$ has been studied intensively in recent years
and is sometimes also referred to as the set covariance.  It follows immediately from basic principles that the covariogram is also 
equal to the {\bf autocorrelation of $1_{A}$}:
\begin{equation}
    g_A(x) := 1_A*1_{-A}(x):=\int_{\R^d} 1_A(t)1_{-A}(x-t)dt,
\end{equation}
where $1_A$ is the indicator function of the body $A$.
More generally, we also define the {\bf cross-correlation} of two bodies $A, B$ as 
$g_{A, B}(x):=1_A*1_{-B}(x):=
\int_{\R^d} 1_A(t)1_{-B}(x-t)dt$, and again an elementary computation gives 
$g_{A, B}(x)=\vol(A\cap (B+n))$.

First, we apply a result of Bombieri \cite{Bombieri} to give universal lower bounds for the series
$\sum_{n \in \Lat} g_A(n)$, and completely characterize the compact sets $A$ that give equality (Corollary \ref{main cor}) in these lower bounds, which turn out to be multi-tiling bodies. 
We recall the classical formula of Siegel, which in itself is also an extension of Minkowski's first theorem (see for example \cite{Robins}).

\begin{thm*}
[C. L. Siegel, 1935]\label{Siegel1}
Let $\PP \subset \R^d$ be a convex body, and $\Lat \subset \R^d$ a full-rank lattice.
If the only lattice point of $\mathcal{L}$ in the interior of 
$\frac{1}{2} \PP-\frac{1}{2} \PP$ is the origin, 
then we have:
\begin{equation}\label{Siegel_lattice}
	2^{d} \operatorname{det} \mathcal{L}=\operatorname{vol} \PP+\frac{4^{d}}{\operatorname{vol} \PP} \sum_{\xi \in \mathcal{L}^{*}\setminus\{0\}}\left|\hat{1}_{\frac{1}{2} \PP}(\xi)\right|^{2}. 
\end{equation}
Taking $\QQ:=\frac{1}{2}\PP$, we may equivalently write:
\begin{equation}
	\operatorname{det} \mathcal{L}=\operatorname{vol} \QQ+\frac{1}{\operatorname{vol} \QQ} \sum_{\xi \in \mathcal{L}^{*}\setminus\{0\}}\left|\hat{1}_{\QQ}(\xi)\right|^{2}. 
\end{equation}
\hfill $\square$
\end{thm*}

\noindent
C. L. Siegel's original proof  of  \eqref{Siegel_lattice} used the Parseval identity \cite{Siegel2}.  Bombieri \cite{Bombieri}
 obtained an extension of Siegel's formula
(again using Parseval's identity)
by relaxing the hypothesis in Siegel's formula \eqref{Siegel_lattice}, so that $\PP$ is allowed to contain exactly $N$ interior lattice points of $\Lat$, with $N>1$.

\medskip
\begin{thm*}[E. Bombieri, 1962]
\label{Bombieri's Theorem}
Let $A\subset \R^d$ 
be a compact set, and $\Lat\subset \R^d$ any full-rank lattice. 
Then:
\begin{equation} \label{Bombieri eq}
   \sum_{n\in\Lat} \vol\big(A\cap\left(A+x+n\right)\big)=
    \frac{1}{\det \Lat} \sum_{\xi \in \Lat^*}
    \left|\hat{1}_{A}(\xi)\right|^{2}
    e^{2\pi i \langle \xi, x \rangle},
\end{equation}
for each $x\in \R^d$.
\hfill $\square$
\end{thm*} 
\noindent
Formula \eqref{Bombieri eq} is called the Bombieri-Siegel formula.  
Here we give an extension of the Bombieri-Siegel formula, to include a more general class of functions, and also to refine the summation index for the series on the left-hand side of \eqref{Bombieri eq}.  As a consequence, we find related inequalities and we classify the cases of equality in terms of multi-tiling bodies.

Our extension of Bombieri's theorem,  namely Theorem \ref{general_case} below, has more general applications to the interaction between any two distinct convex bodies. 
We begin with the following known variation of the Poisson summation formula (PSF).  This particular variation of Poisson summation doesn't seem to as well-known as some of the other variations, so we include an independent proof in Appendix \ref{Known lemmas}.  
\begin{thm*}[Poisson summation \cite{RichardStrungaru}] 
Suppose that $g:\R^d \rightarrow \C$ is compactly supported, continuous, and  
$\hat g \in L^1(\R^d)$.
Then we have:
\begin{equation} 
\label{OurPoissonSummation}
      \sum_{n\in\Lat}g(n+x)
     =\frac{1}{\det{\Lat}}\sum_{m\in\Lat^*}\widehat{g}(m)e^{2\pi i \langle m, x\rangle},
\end{equation}
for any full-rank lattice $\Lat$, and all $x \in \R^d$.
Here the equality holds pointwise, 
and both series converge absolutely and uniformly to continuous functions.
\hfill 
\hyperlink{Proof of Our Poisson Summation}{\rm{(Proof)}}
$\square$
\end{thm*}
 \noindent
If identity \eqref{OurPoissonSummation}
holds for $g(x)$,
then we'll call $g$  
{\bf Poisson summation friendly}.
We include an independent proof of the latter Poisson summation formula (see Appendix), because it is not as 
well-known as most of the other variants of PSF.  Given any measurable sets $A, B\subset \R^d$, if the function 
\begin{equation}
    g(x):=\left(1_{A} * 1_{B}\right)(x),
\end{equation}
is Poisson summation friendly, then we say that the {\bf sets $A$ and $B$ are 
 Poisson summation friendly}. 
It is very natural to wonder how general the class of such pairs of PSF sets can be, because they arise naturally in the study of  cross covariograms. Our first main result is the following slight extension of Bombieri's Theorem, namely  Theorem \ref{general_case}, part 
\eqref{third part of main result}.  

\medskip
\begin{thm}\label{general_case}
We fix any two measurable, bounded sets 
$A, B\subset \R^d$, and we let $f, g \in L^1(\R^d)\cap L^2(\R^d)$ enjoy the  following properties:
both functions are bounded, $f$ is compactly supported on $A$, and 
$g$ is compactly supported on $B$. 
For any full-rank lattice $\Lat\subset \R^d$ we have:
\begin{enumerate}[(a)]
\item  
\label{first part of main result}
$f*g$  is  continuous on $\R^d$, and is also Poisson summation friendly.
\item 
\label{second part of main result}
Consequently:
\begin{equation} 
\label{generalized convolution identity}
   \sum_{n\in\Lat} 
   (f*g)(x+n) 
=\frac{1}{\det \Lat} \sum_{\xi \in \Lat^*}
   \hat{f}(\xi)
 \hat{g}(\xi) 
    e^{2\pi i \langle \xi, x \rangle},
\end{equation}
for each $x\in \R^d$.
\smallskip
\item 
\label{third part of main result}
If we choose $f := 1_A$ and 
$g:= 1_{-B}$, then as a special case of part \eqref{second part of main result}, we get:
\begin{equation} 
\label{extending eq. of Bombieri}
\sum_{n\in\Lat}
\vol\Big(
A\cap\left(B+x+n\right)
\Big)=
\frac{1}{\det \Lat} \sum_{\xi \in \Lat^*}
   \hat{1}_{A}(\xi)
 \overline{ \hat{1}_{B}(\xi) }
    e^{2\pi i \langle \xi, x \rangle},
\end{equation}
for each $x\in \R^d$.
\hfill 
\hyperlink{Proof of general_case}{\rm{(Proof)}}
$\square$
\end{enumerate}
\end{thm} 

In particular, Theorem \ref{general_case}, 
part \eqref{first part of main result} 
tells us that {\bf every bounded, measurable set 
$B\subset\R^d$  is a Poisson summation friendly set}, because in this case the function $f(x):= 1_B*1_{-B}$ is a Poisson summation friendly function.

The Bombieri-Siegel formula \eqref{Bombieri eq} is the case 
$f= g$, and hence $A=B$.
A trivial observation is that 
if a centrally symmetric body $\BB$ contains exactly $N$ interior lattice points then $N$ must be odd.  The reason is easy:
for all nonzero $n\in\operatorname{int}(\BB)\cap\Lat$, we have  $-n\in\operatorname{int}(\BB)\cap\Lat$. Including the origin, we therefore have an odd number of lattice points. 

But this brings up another question: can we describe more precisely  the lattice sum on the left of
\eqref{Bombieri eq} ? We answer this question in the affirmative, when working with compact sets.  Namely,
we obtain another  refinement of Bombieri's identity \eqref{Bombieri eq}, by proving that for compact sets the sum on the 
left-hand side of  
\eqref{Bombieri eq} may be restricted to precisely the lattice points contained in the interior of $\QQ - \QQ$. More generally, we may restrict the lattice sum
on the left of   \eqref{Bombieri eq} 
to the interior of $A-B$, for any pair of compact sets 
$A, B$, as the following theorem shows.  

Our second main result is the following extension of Bombieri's Theorem, which gives a precise finite index of summation. Its proof uses the technical 
Lemma~\ref{first lemma} below, and due to this Lemma we will assume henceforth that all of our sets are compact, so that in particular they are also Poisson summation friendly.

\begin{thm}[Refined Bombieri-Siegel formula] 
\label{main, for Q}
Let $A, B$ be compact sets in $\R^d$, and let
$\Lat\subset \R^d$ be a full-rank lattice.
\begin{enumerate}[(a)]
\item  \label{part a of Theorem 3} 
Then we have
\begin{equation}
   \sum_{n \in \Int(A-B) \cap \Lat } \vol\Big(A\cap\left(B+n\right)\Big)
    =\frac{1}{\det \Lat} 
    \vol A \vol B
    + \frac{1}{\det \Lat}
    \sum_{\xi \in \Lat^* \setminus \{0\}}
    \hat{1}_{A}(\xi)
    \overline{\hat{1}_{B}(\xi)},
\end{equation}
\item \label{part b of Theorem 3}
If $\PP$ is centrally symmetric and convex, then:
\begin{equation}\label{Main2}
    \sum_{n \in \Int(\PP) \cap \Lat } \vol\Big(\PP\cap\left(\PP+2n\right)\Big)
    =\frac{1}{\det \Lat}  \frac{\vol^{2} \PP}{2^{d}}
    + \frac{2^d}{\det \Lat}
    \sum_{\xi \in 
    \Lat^* \setminus\{0\}}\left|\hat{1}_{\frac{1}{2} \PP}(\xi)\right|^{2}.
\end{equation}
\end{enumerate}
\hfill 
\hyperlink{Proof of main, for Q}{\rm{(Proof)}}
$\square$
\end{thm}
\bigskip
We say that a body $\PP\subset \R^d$ {\bf $k$-tiles by translations with a lattice} $\Lat \subset \R^d$
if
\begin{equation}
    \sum_{n\in \Lat} 1_{\PP+n}(x) = k,
\end{equation}
for all $x \in \R^d$, except for 
$x\in \partial \PP+ \Lat$, a measure zero set. When $k=1$, this is the classical definition of tiling by translations, where there are no overlaps between the translated interiors of $\PP$.  But when $k>1$, we note that for such a $k$-tiling  the translates of $\PP$ will always overlap each other.  The field of $k$-tiling has recently experienced a renaissance, and one of the earliest works relating $k$-tiling to Fourier analysis was done by M. Kolountzakis \cite{Kolountzakis1}.

Next, we discretize our results above, and consider any 
{\em finite} set of integer points $F\subset \Z^d$, of cardinality $|F|$.  Setting $\square:= [-\tfrac{\varepsilon}{2}, \tfrac{\varepsilon}{2}]^d$,
we now ``thicken'' each point $n\in F$, by replacing it with a cube of sidelength 
$\varepsilon$ centered at $n$. In other words, we replace each integer point $n \in F$ by the little $d$-dimensional cube
$ n+[-\tfrac{\varepsilon}{2}, \tfrac{\varepsilon}{2}]^d:= n+\square$, for a fixed $0 < \varepsilon \leq 1$. 
Hence we've modified the finite set $F$ into the $d$-dimensional compact set 
\begin{equation}
{A}:=  \bigcup_{n\in F} 
\left( \square+n \right),
\end{equation}
which we call the {\bf $\varepsilon$-thickening} of $F$. 
Moreover, for $\varepsilon \in (0, 1]$, the $\varepsilon$-cubes are centered at integer points and therefore do not overlap in any set of positive measure,
allowing us to conclude 
that $\vol A = \varepsilon^d |F|$.

\begin{thm}[Discretized Bombieri-Siegel formula] \label{application: arithmetic combintorics 1}
Let $F  \subset \Z^d$ 
be any finite set of integer points, and fix a full-rank integer sublattice $\Lat\subset \Z^d$.
Let $A$ be the $\varepsilon$-thickening of $F$,
for each fixed $\varepsilon \in (0, 1]$. 
\begin{enumerate}[(a)]
\item  \label{first eq. of AC proof}
We have:
\begin{equation}
   \sum_{n\in (F-F)\cap \Lat} 
   |F\cap (F+n)| 
  =\frac{1}{\det \Lat} |F|^2\varepsilon^{d} 
  +
  \frac{1}{\varepsilon^d \det \Lat} 
  \sum_{\xi \in \Lat^*\setminus 0}
    \left | \hat{1}_{A}(\xi)\right |^2.
\end{equation}
\item  \label{second eq. of AC proof}   
Equivalently, we have the explicit respresentation:
\begin{equation}  
      \sum_{n\in (F-F)\cap \Lat} 
   |F\cap (F+n)| 
  =\frac{1}{\det \Lat} |F|^2\varepsilon^{d}
  +
 \frac{1}{ \det \Lat} 
  \sum_{\xi \in \Lat^*\setminus 0} 
\left(
\prod_{k=1}^d \sinc^2(\pi \varepsilon \xi_k)
 \left|
 \sum_{n \in F}  e^{2\pi i\langle  \xi, n\rangle}
 \right|^2
 \right),
 \end{equation}
 where 
 $\sinc(\pi t):=
 \begin{cases}
  \frac{\sin(\pi t)}{\pi t} & \text{ if } t \not=0, \\
  1 & \text{ if } t =0.
 \end{cases}$
\end{enumerate}
\hfill  
\hyperlink{first theorem for integer sets}{\rm{(Proof)}} $\square$
\end{thm}



In Theorem \ref{tiling the integer lattice with a finite set of integers}
 below, we will be interested in letting $\varepsilon =1$.  
Next, we apply Theorem \ref{main, for Q}, 
part  \eqref{part a of Theorem 3}, with $A=B$, where $A$ is by definition the $1$-thickening of a finite set of integer points 
$F\subset \Z^d$:
\begin{equation}\label{special discrete case 1}
   \sum_{n \in \Int(A-A) \cap \Lat } \vol\Big(A\cap\left(A+n\right)\Big)
    = 
    \frac{(\vol A)^2}{\det \Lat}
    +
    \frac{1}{\det \Lat}
    \sum_{\xi \in  
    \Lat^* \setminus\{0\}}
    \left | \hat{1}_{A}(\xi)\right |^2.
\end{equation}
We recall the easy fact that 
$ F\cap (F+n) \not= \varnothing
\iff n \in F - F$. By assumption we have $F\subset \Z^d$, so that 
$F-F\subset \Z^d$ as well.

We now fix a sublattice $\Lat\subset \Z^d$.
From the discussion above, we see that for each   
$n \in \Lat\cap(F-F)$, we have the following relation between $F$ and its 
$\varepsilon$-thickening $A$:
\begin{equation}
\vol\Big(A\cap\left(A+n\right)\Big)
= |F\cap (F+n)| \varepsilon^d,
\end{equation}
because any integer translate of $A$ intersects $A$ itself in an integer number of translated copies of the 
$\varepsilon$-cube
$\square:=\left[-\tfrac{\varepsilon}{2}, 
\tfrac{\varepsilon}{2}\right]^d$.  
So we can rephrase equation  \eqref{special discrete case 1} 
in terms of the finite set $F\subset \Z^d$.

\begin{cor}
\label{lower bound for A-A}
Let $F \subset \Z^d$ 
be any finite set of integer points, and
fix $\Lat \subset \Z^d$, a full-rank sublattice.  
Then 
\begin{equation} \label{lower bound for finite covariograms}
   \sum_{n\in (F-F)\cap \Lat} 
   |F\cap (F+n)| 
  \geq
  \frac{1}{\det \Lat} |F|^2.
\end{equation}
\end{cor}
\begin{proof}
The result follows directly from 
Corollary \ref{application: arithmetic combintorics 1}, 
part \eqref{first eq. of AC proof}, with $\varepsilon =1$.
\end{proof}

Theorem \ref{lower bound for A-A}  is number-theoretic in the sense that here $\det \Lat$ and $|F|$ are positive integers.
A natural question arises: what are the cases of equality in \eqref{lower bound for finite covariograms}?
There is a particularly interesting answer, which we give in Theorem 
\ref{tiling the integer lattice with a finite set of integers} below. 
We'll make use of the following trivial but useful observation.
\begin{observation*}\label{observation for tiling}
A finite collection of integer points $F\subset\R^d$ multi-tiles $\Z^d$ if and only if  its
$1$-thickening  $A$ multi-tiles $\R^d$ with the same multiplicity.
\end{observation*} 

\begin{thm}
\label{tiling the integer lattice with a finite set of integers}

Let $F \subset \Z^d$ 
be any finite set of integer points, and
fix $\Lat \subset \Z^d$, a full-rank sublattice. Then the following statements are equivalent. 
\begin{enumerate}[(a)] 
   \item \label{part a of lattice tilings}
  The finite set $F$ multi-tiles $\Z^d$ by  translations with $\Lat$, and with multiplicity $k=\frac{|F|}{\det\Lat}$.
  \item  \label{part b of lattice tilings}
\begin{equation}
\label{tiling Z^d}
 \sum_{n\in (F-F)\cap \Lat} |F\cap (F+n)| 
  = k|F|,
\end{equation}
with $k=\frac{|F|}{\det\Lat}$.
\item \label{part c of lattice tilings} For each nonzero $\xi \in \Lat^*$
such that none of the coordinates of $\xi$ are integers, the following exponential sum vanishes:
\begin{equation}
     \sum_{n \in F}  
     e^{2\pi i\langle  \xi, n\rangle}=0.
\end{equation}
\end{enumerate}
\hfill  
\hyperlink{main theorem for integer lattices}{\rm{(Proof)}} $\square$
\end{thm}
Returning to the continuous context, we may now give 
several corollaries of Theorem \ref{main, for Q},
the refined Bombieri-Siegel formula. 
The following inequality gives us a best-possible lower bound for the sum of the covariogram over a lattice.  Interestingly,  the lower bound in equation \eqref{main inequality} 
below is achieved by a body $\QQ$ precisely when $\QQ$ is a $k$-tiling polytope. 

\begin{cor} \label{main cor}
Let $A \subset \R^d$ be a compact set and let 
$\Lat\subset \R^d$ be a full-rank lattice.   Then we have:
\begin{equation} \label{main inequality}
\sum_{n\in\Lat\cap \interior(A-A)} \vol\Big(A\cap\left(A+n\right)\Big) 
    \geq  \frac{\vol^2 A}{\det\Lat}.
\end{equation}
\hfill  
\hyperlink{proof of main cor}{\rm{(Proof)}} $\square$
\end{cor}


\begin{cor} \label{multi-tiling Cor}
Let $A \subset \R^d$ be a compact set and let 
$\Lat\subset \R^d$ be a full-rank lattice.   Then the following are equivalent:
\begin{enumerate}[(a)]
    \item 
    \begin{equation} 
    \sum_{n\in\Lat\cap \interior(A-A)} \vol\Big(A\cap\left(A+n\right)\Big) 
   = \frac{\vol^2 A}{\det\Lat}.
\end{equation}
    \item 
    The body $A$ $k$-tiles $\R^d$ with the lattice $\Lat$, where 
    $k=\frac{\vol A}{\det\Lat}$.
    \item 
\bigskip
For all compact sets $B\subset \R^d$, 
we have
\begin{equation}\label{part (c) of Cor 2}
     \sum_{n\in\Lat\cap \interior(A-B)} \vol\Big(A\cap\left(B+n\right)\Big) 
   = \frac{1}{\det \Lat} \vol A \vol B.
\end{equation}
\end{enumerate}
\hfill 
\hyperlink{proof of multi-tiling Cor}{\rm{(Proof)}}
$\square$
\end{cor}


Another way to conceptualize Corollary \ref{main cor} is by recalling a well-known identity, as follows. 
With $f(x):= (1_{\QQ}*1_{-\QQ})(x) = \vol\Big(\QQ \cap (\QQ + x)\Big)$, we have $\hat f(0) = | \hat 1_\QQ(0)|^2 =  \vol^2(\QQ)$.  But we also have, via the inverse Fourier transform,  $\hat f(0) = \int_{\R^d} f(x) dx$, giving us the known identity
\begin{equation}
  \int_{\R^d} \vol\Big(Q\cap(Q+x)\Big) dx=\vol^2 \QQ. 
\end{equation}
Putting this together with the right-hand side of equation 
\eqref{main inequality}, we've just proved the following reformulation of Corollary \ref{main cor}, which appears to resemble a quadrature formula. 

\begin{cor} \label{main cor reformulated}
Let $\QQ \subset \R^d$ be a compact set and let $\Lat\subset \R^d$ be a full-rank lattice.
Then we have: 
\begin{equation} \label{secondary main inequality}
    \sum_{n\in\Lat\cap \interior(\QQ-\QQ)} \vol\Big(\QQ\cap(\QQ+n)\Big)
    \geq  \frac{1}{\det\Lat}
    \int_{\R^d} \vol\Big(\QQ\cap(\QQ+x)\Big) dx.
\end{equation}
Moreover, equality occurs in 
\eqref{secondary main inequality} $\iff$
$\QQ$ $k$-tiles $\R^d$ with the lattice $\Lat$.
\hfill $\square$
\end{cor}

\bigskip
We also obtain the following interesting spectral identity for the product of two volumes of any compact sets, provided we pick a sufficiently sparse lattice $\Lat$.

\begin{cor} [Spectral formula for the product of two volumes]
\label{max of FT is at the origin}
Let $\Lat\subset \R^d$ be a full-rank lattice, and let $A, B \subset \R^d$ be compact sets possessing the following property. 
We assume that $A$ is disjoint from the nonzero lattice translates $\{B+n \mid n \in\Lat\}$.
In addition, let $x\in \R^d$ be a nonzero vector with the property that both $A$ and $B+x$ are also mutually disjoint from their nonzero translations by  vectors from $\Lat$.  
Then 
\begin{equation}
\vol\left(A\right) \vol\left(B\right) 
=
 -\sum_{\xi \in \Lat^* \setminus \{0\}}
   \hat{1}_{A}(\xi)
   \overline{\hat 1_{B}(\xi)}
    \cos(2\pi \langle \xi, x \rangle).
\end{equation}
\hfill 
\hyperlink{proof of max of FT is at the origin}{\rm{(Proof)}}
$\square$
\end{cor}

\noindent
The special case $A=B$ of Corollary 
\ref{max of FT is at the origin} is  interesting in its own right:
\begin{equation}
    \vol^2 A
=
 -\sum_{\xi \in \Lat^* \setminus \{0\}}
   \left| \hat{1}_{A}(\xi) \right|^2
    \cos(2\pi \langle \xi, x \rangle),
\end{equation}
and in particular the latter identity 
holds for any compact set $A$ that satisfies the hypothesis of Corollary \ref{max of FT is at the origin}.  We remark that given any compact sets $A, B \subset \R^d$, it's quite easy to find a lattice $\Lat$ that satisfies the hypotheses of Corollary \ref{max of FT is at the origin}.

We also give applications to arithmetic combinatorics, in Section \ref{sec:arithmetic combinatorics}, by discretizing the Bombieri-Siegel formula and its extension above, to handle any finite set of integer points in $\Z^d$.  In particular, Corollaries \ref{application: arithmetic combintorics 1}, 
 \ref{lower bound for A-A},  
 \ref{tiling the integer lattice with a finite set of integers} and 
  \ref{discrete Minkowski} give us the relevant discretizations.

A related research direction is the question of using a finite set of integer points $F$ to multi-tile $\Z^d$ by any possible set  $S$ of integer translations.  The {\bf periodicity conjecture} states that it is always sufficient to let $S$ be  lattice.  Bhattacharya \cite{Bhattacharya}  recently proved the periodicity conjecture for dimension $2$.  However,
there is a recent breakthrough on the periodicity conjecture by Greenfeld and Tao \cite{GreenfeldTao}, disproving it for all sufficiently large dimensions $d$.

Another strong motivation for studying the covariogram is the following conjectural characterization of a convex body $K$, posed by G. Matheron
\cite{Matheron}, and studied intensively over the past $30$ years by G. Bianchi, as well as other researchers \cite{Bianchi1} \cite{Bianchi2}.

\begin{conj}[Matheron, 1986]
\label{conj:Matheron}
The covariogram $g_K$ determines a convex body $K$, among all convex bodies, up to translations and reflections.
\hfill $\square$
\end{conj} 

We note that knowledge of $\hat 1_K(\xi)$, for all $\xi \in \R^d$, uniquely determines $K$, for any convex body $K$; in other words, the Fourier transform is a complete invariant for any body $K$.   Hence Conjecture \ref{conj:Matheron} is equivalent to saying that knowledge of $|\hat 1_K|$ determines $K$ (up to translations and reflections), because  $\hat g_K := \widehat{1_K*1_{-K}} =|\hat 1_K|^2$. 

It is known, for example,  that centrally symmetric convex bodies in any dimension are determined by their covariogram, up to translations \cite{Bianchi2}. 

\bigskip
\section{Preliminaries}
The following somewhat technical result is required for the proof of Theorem \ref{main, for Q}.
\begin{lem}\label{first lemma}
Let $A, B\in \R^d$ be compact sets, 
and let 
$\Lat\subset \R^d$ be a full-rank lattice.  
\begin{enumerate}[(a)]
    \item 
For any 
$n\in\Lat\cap\partial(A-B)$,
we have 
\begin{equation}
    \vol{\Big(A\cap\left(B+n\right)\Big)}=0.
\end{equation}
\item \label{part b of first lemma} Consequently, we also have:
\begin{equation}
 \sum_{n\in\Lat} \vol\Big(A\cap\left(B+n\right)\Big)
 =  
 \sum_{n\in\Lat\cap \interior(A-B)} \vol\Big(A\cap\left(B+n\right)\Big).
\end{equation}
\end{enumerate}
\hfill 
\hyperlink{Proof of first lemma}{\rm{(Proof)}}
$\square$
\end{lem}
We define 
 the function spaces $L^p(\R^d):=\{f:\R^d \rightarrow \C \mid 
\int_{\R^d} |f(x)|^p dx < \infty\}$, for each $1\leq p < \infty$.  Here we'll use $p=1$ and $p=2$.  Our Fourier transforms are defined by the traditional convention: 
$\hat f(\xi):= \int_{\R^d}f(x) 
e^{-2\pi i \langle x, \xi \rangle} dx$, when the integral converges.  A full-rank lattice $\Lat \subset \R^d$ is a lattice that has a $d$-dimensional basis.
We begin by mentioning an elementary argument that proves the following fact.

\begin{lem} \label{simple lemma} If a body $\QQ$ $k$-tiles $\R^d$ we have the identity
\begin{equation}
    \sum_{n \in\Lat} \vol\left(\QQ\cap\left(\QQ+n+x\right)\right)=k\cdot\vol\QQ
\end{equation}
for every fixed $x\in\R^d$.
\end{lem}
\begin{proof}
For a $k$-tiling we must have
\begin{equation}
    \sum_{n\in \Lat}1_{\QQ+n}(y)=k
\end{equation}
for all $y\in\R^d$, except those points $y$ that lie on the boundary of $\QQ$ or its translates under $\Lat$.
Therefore
\begin{align*}
    \sum_{n \in\Lat} \vol\left(\QQ\cap\left(\QQ+n+x\right)\right)
    &=\sum_{n \in \Lat} \int_{\mathbb{R}^{d}} 1_{\QQ\cap\left(\QQ+n+x\right)}(y)  d y\\
    &=\sum_{n \in \Lat} \int_{\mathbb{R}^{d}} 1_{\QQ}(y)\cdot1_{\QQ+n+x}(y)  d y\\
    &=\sum_{n \in \Lat} \int_{\mathbb{R}^{d}} 1_{\QQ-x}(y)\cdot1_{\QQ+n}(y)  d y\\
    &=\int_{\mathbb{R}^{d}} 1_{\QQ-x}(y)\cdot\sum_{n \in \Lat}1_{\QQ+n}(y)  d y\\
    &=k\cdot\vol\QQ,
\end{align*}
where we used the translation-invariance of the integral and exchanged the sum with the integral since the integral is nonzero only for a finite number of $n \in \Lat$.
\end{proof}

It turns out that the converse of Lemma \ref{simple lemma} is also true, and is part of
of Corollary \ref{multi-tiling Cor}.    
A very useful Fourier equivalence for translational tilings is given by the following result, 
due to M. Kolountzakis 
\cite{Kolountzakis1} (see also \cite{Robins}, Theorem 5.5). 

\medskip
\begin{thm}[Kolountzakis]\label{Kol}

  Let $\PP$ be a body in $\R^d$, and let $\Lat\subset \R^d$ be a full-rank lattice.  
  Then the following conditions are equivalent:
\begin{enumerate}
    \item  \label{Kol, part 1}
    $\hat{1}_{\PP}(\xi)=0$, for all nonzero $\xi\in\Lat^*$, the dual lattice.
    \item    \label{Kol, part 2}
    $\PP$ $k$-tiles $\R^d$ by translations with $\Lat$.
    \end{enumerate}
Either of the above conditions already 
implies that
    $k=\frac{\vol{P}}{\det{\Lat}}$,  a positive integer.
\hfill $\square$
\end{thm}

\begin{rem}
We may rewrite equation 
\eqref{Bombieri eq} as
\begin{equation} \label{real and imaginary}
   \det \Lat\sum_{n\in\Lat} \vol\left(\QQ\cap\left(\QQ+n+ x\right)\right)=
     \sum_{\xi \in \Lat^*}
    \left|\hat{1}_{\QQ}(\xi)\right|^{2}
    \cos{(2\pi\langle \xi, x \rangle)}+i\sum_{\xi \in \Lat^*}
    \left|\hat{1}_{\QQ}(\xi)\right|^{2}\sin{(2\pi\langle \xi, x \rangle)},
\end{equation}
from which we may conclude that the imaginary part vanishes:
\begin{equation}\label{imaginary part is zero}
    \sum_{\xi \in \Lat^*}
    \left|\hat{1}_{\QQ}(\xi)\right|^{2}\sin{(2\pi\langle \xi, x \rangle)}=0,
\end{equation}
for all $x\in\R^d$. It's also easy to see that this equality \eqref{imaginary part is zero} 
 follows from first principles.  Namely, we first observe that 
$ \left|\hat{1}_{\QQ}(\xi)\right|^{2} 
= \hat{1}_{\QQ}(\xi) \overline{\hat{1}_{\QQ}(\xi)} 
= \hat{1}_{\QQ}(\xi) \hat{1}_{\QQ}(-\xi)$.  Now,  
we may split $\R^d$ into $3$ regions defined by the hyperplane orthogonal to $x$, namely $H:= \{ \xi \in \R^d \mid \langle x, \xi \rangle = 0\}$,
and the two half-spaces defined by $H$. 
Each $\xi \in \Lat^*$ on one side of $H$ gives a summand that cancels with the corresponding $-\xi$ on the other side of $H$, while the $\xi$'s that lie on $H$ yield a vanishing summand. 

The above discussion shows that 
\begin{equation}\label{alternate main Theorem}
   \sum_{n\in\Lat} \vol\left(\QQ\cap\left(\QQ+n+ x\right)\right)=
    \frac{1}{\det \Lat}  
    \sum_{\xi \in \Lat^*}
    \left|\hat{1}_{\QQ}(\xi)\right|^{2}
    \cos{(2\pi\langle \xi, x \rangle)},
\end{equation}
valid for each $x\in \R^d$, is equivalent to equation
\eqref{Bombieri eq} 
of Theorem  
\ref{Bombieri's Theorem}, a form which Bombieri already 
had in \cite{Bombieri}.
\hfill $\square$
\end{rem}


\section{Examples of the discretized Bombieri-Siegel formula}
\label{sec:arithmetic combinatorics}

\begin{example}
We let $F:= \left\{ (0,0), (2, 1), (1, 2), (1, 3) \right\}$, and let $\Lat$ be the lattice defined by the basis 
 $\{ (1, 1), (-2, 2)\}$.  
It can be verified that $F$ tiles $\Z^2$ by translations with $\Lat$. 
\hfill $\square$
\end{example}

As an immediate consequence of Corollary 
\ref{lower bound for A-A}
and Theorem \ref{tiling the integer lattice with a finite set of integers}, we next obtain a discrete version of Minkowski's formula, as follows. 
\begin{cor}
\label{discrete Minkowski}    
Let $F \subset \Z^d$ 
be any finite nonempty set of integer points, and
fix $\Lat \subset \Z^d$, a full-rank sublattice.  
If the only lattice point of $\Lat$ in
$F-F$ is the origin, then 
\begin{equation} \label{Discrete Minkowski identity}
\det\Lat\geq|F|.
\end{equation}
Moreover, equality occurs in 
\eqref{Discrete Minkowski identity}
if and only if the finite set 
$F$ multi-tiles
$\Z^d$ by translations with $\Lat$, and with multiplicity $k=\frac{|F|}{\det\Lat}$.
\end{cor}
\begin{proof}
The inequality of \eqref{lower bound for finite covariograms} yields 
\begin{equation}
   \sum_{n\in (F-F)\cap \Lat} 
   |F\cap (F+n)|=\sum_{n=0} 
   |F\cap (F+n)| =|F|\geq \frac{1}{\det \Lat} |F|^2,
\end{equation}
which gives us $\det\Lat\geq|F|$.  The equality case follows immediately from Theorem \ref{tiling the integer lattice with a finite set of integers}.
\end{proof}
It is worth noting that Corollary \ref{discrete Minkowski} can also be proved in an elementary way by reducing the points of $F$ modulo a fundamental domain of the lattice $\Lat$.

\begin{example}
\rm{
Consider the set $F:= \{ 1, 3, 4, 6\}$, which happens to tile $\Z$ by the set of translations that belong to the lattice $\Lat:=4\Z$, with multiplicity $k=1$ (using brute-force).  
Here 
\[
F-F = \{ -5, -3, -2, -1, 0, 1, 2, 3, 5\}.
\]
To confirm the latter tiling claim by using Theorem \ref{tiling the integer lattice with a finite set of integers}, we compute 
\begin{align*}
\sum_{n\in (F-F)\cap 4\Z} |F\cap (F+n)|
&=|F\cap (F+0)|= 4,
\end{align*}
whereas $k|F| = 4$ as well, confirming that according to Theorem \ref{tiling the integer lattice with a finite set of integers}, the finite set $F$ of integers must tile $\Z$, with multiplicity $k=1$.
}
\hfill $\square$
\end{example}

\begin{example}
\rm{
Consider the set of points in $\Z^2$:
\begin{equation*}
 F:=\{(0, 0), (0, 2), (-1, 3), (1, 3), (-1, -3), (0, -2), (1, -3), (1, 
  1), (-1, -1)\}   
\end{equation*}
and the integer sublattice $\Lat$ defined by the basis  $\left\{\tvect{3}{0},\tvect{1}{1}\right\}$, as shown in Figure \ref{example_discrete_case}.   

\begin{figure}[ht]
    \centering
\includegraphics[width=7cm]{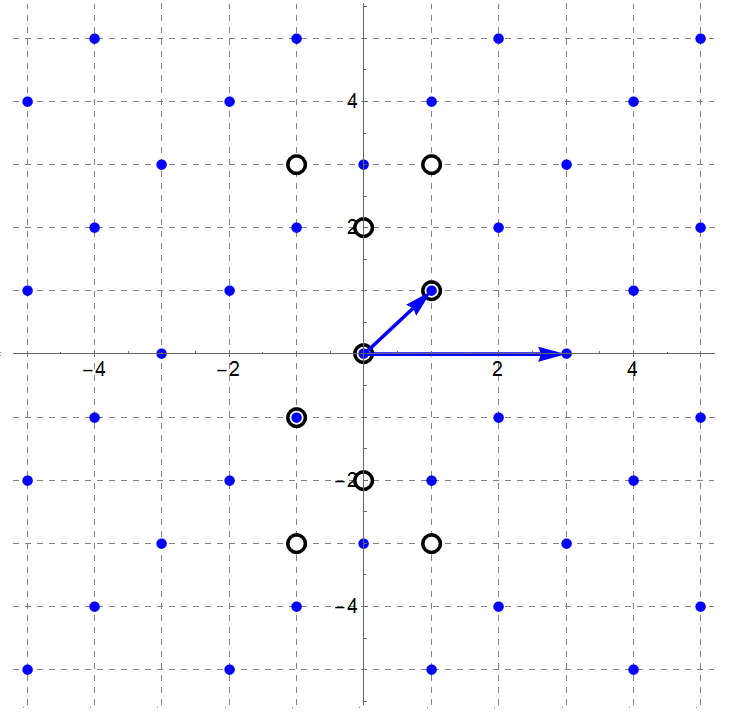}
    \caption{The set $F$ is shown by the Black circles. The integer sublattice $\Lat$ is shown by the blue dots. A basis for $\Lat$ is shown by the two arrows.} \label{example_discrete_case}
\end{figure}
\noindent
Here the difference set $F-F$ is given by:
\begin{align*}
    F-F =& 
    \Big\{
    (-2, -6), (-2, -4), (-2, -2), (-2, 0), (-2, 2), (-2, 
  6), (-1, -5), (-1, -3), (-1, -1),\\
  & (-1, 1), (-1, 3), (-1, 
  5), (0, -6), (0, -4), (0, -2), (0, 0), (0, 2), (0, 4), (0, 
  6),(1, -5),\\
  &(1, -3), (1, -1), (1, 1), (1, 3), (1, 5), (2, -6), (2, -2), (2, 0), (2, 2), (2, 4), (2, 6)
  \Big\}.
\end{align*}
We therefore have:
\begin{equation}\label{nj}
    (F-F)\cap \Lat = \{(-2, -2), (-1, -1), (-1, 5), (0, -6), (0, 0), (0, 6), (1, -5), (1, 1), (2, 2)\},
\end{equation}
which is the index of summation for the left-hand-side of identity
\eqref{tiling Z^d}, as in Figure \ref{example_discrete_case_2}.
Labelling the latter elements by
$n_j\in (F-F)\cap\Lat$, we have:
\begin{align*}
    F\cap (F+n_1)= &\{(-1,-1)\};\\
    F\cap (F+n_2)= &\{(-1,-3),(-1,-1),(0,0),(0,2)\};\\
    F\cap (F+n_3)= &\{(-1,3),(0,2)\};\\
    F\cap (F+n_4)= &\{(-1,-3),(1,-3)\};\\
    F\cap (F+n_5)= &\{(-1,-3),(-1,-1),(-1,3),(0,-2),(0,0),(0,2),(1,-3),(1,1),(1,3)\};\\
    F\cap (F+n_6)= &\{(-1,3),(1,3)\};\\
    F\cap (F+n_7)= &\{(0,-2),(1,-3)\};\\
    F\cap (F+n_8)= &\{(0,-2),(0,0),(1,1),(1,3)\};\\
    F\cap (F+n_9)= &\{(1,1)\}.
\end{align*}
\begin{figure}[ht]
    \centering
    \includegraphics[width=8cm]{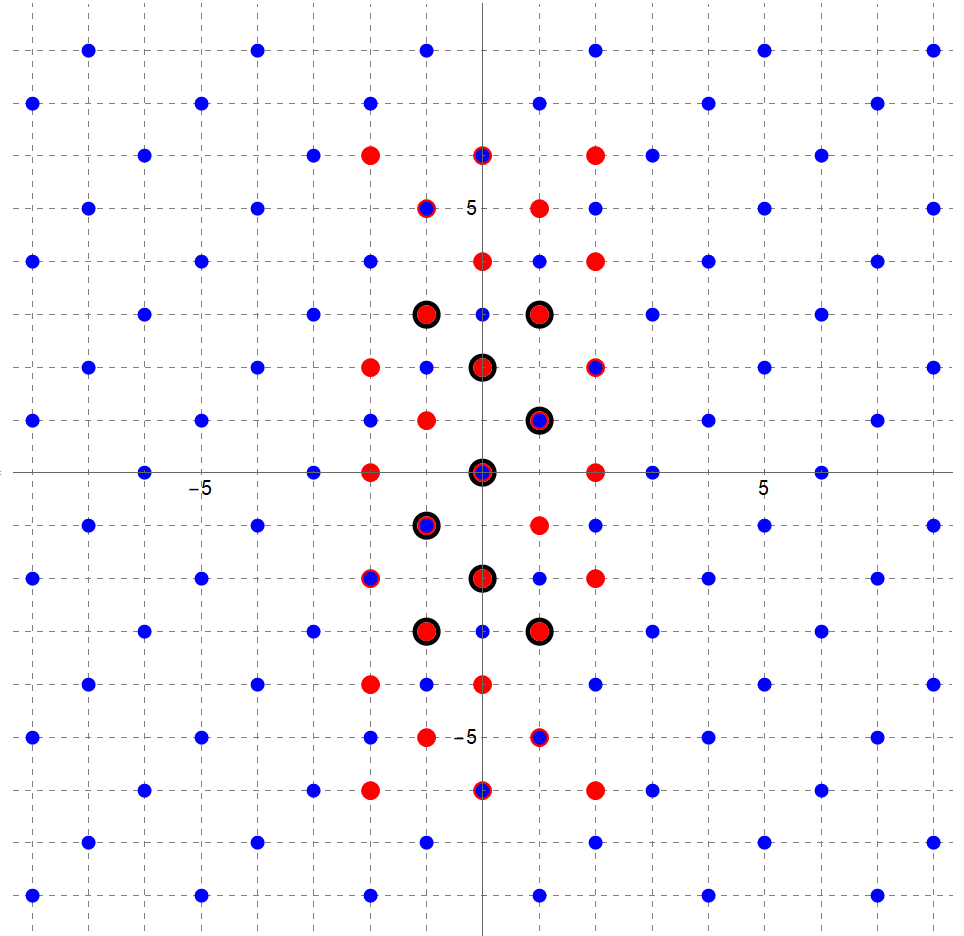}
    \caption{The set $F$ is shown by the Black circles. The integer sublattice $\Lat$ is shown by the blue dots. The difference set $F-F$ is displayed by the $31$ red dots.} \label{example_discrete_case_2}
\end{figure}
Finally, we compute the right-hand side of \eqref{tiling Z^d}:
\begin{align*}
    \sum_{n\in (F-F)\cap \Lat} |F\cap (F+n)|& =  |F\cap (F+n_1)|+\cdots+|F\cap (F+n_9)|\\
    & = 1+4+2+2+9+2+2+4+1 =27.
\end{align*}
On the other hand, the definition 
\begin{equation*}
    k=\frac{|F|}{\det\Lat} = \frac{9}{3\cdot1}=3
\end{equation*} 
gives us the equality
\begin{equation*}
    \sum_{n\in (F-F)\cap \Lat} |F\cap (F+n)|= \frac{|F|}{\det\Lat}\cdot |F| = 27.
\end{equation*} 
Hence Corollary \eqref{tiling the integer lattice with a finite set of integers} part \eqref{part b of lattice tilings} is satisfied for $F$ and the lattice $\Lat$,
implying that the finite set $F$ multi-tiles $\Z^2$ with multiplicity $k=3$.
}\hfill $\square$
\end{example}

\section{A special case: Van der Corput's inequality}

In this short section, we first show that Van der Corput's classical inequality for compact bodies is an immediate consequence of  Corollary~\ref{main cor}. 

\begin{thm}[Van der Corput, 1935]
For any body $\QQ$ and full-rank lattice $\Lat \subset \R^d$, we have
\begin{enumerate}[(a)]
\item \label{first part of Van der Corput}
\begin{equation}\label{Van der Corput inequality}
\#\left\{
     \Int(\QQ-\QQ) \cap \Lat
    \right\} 
    \geq 
    \frac{ \vol \QQ }{ \det \Lat}.
\end{equation}
\item \label{second part of Van der Corput}
If $\QQ$ is also assumed to be convex and centrally symmetric, then we have the inequality 
\begin{equation}\label{easier Van der Corput}
\#\left\{
     \Int(\QQ) \cap \Lat
    \right\}
\geq
\frac{1}{2^d}
\frac{ \vol \QQ }{ \det \Lat}.
\end{equation}
\end{enumerate}
\hfill $\square$
\end{thm}
To see how  inequality \eqref{Van der Corput inequality} of part \eqref{first part of Van der Corput} follows from Corollary~\ref{main cor}, we just note that
\begin{align}
   \#\left\{
     \Int(\QQ-\QQ) \cap \Lat
    \right\}    \vol \QQ 
&=
  \sum_{n \in \Int(\QQ-\QQ) \cap \Lat } \vol\QQ \\
&\geq 
\sum_{n \in \Int(\QQ-\QQ) \cap \Lat } \vol\left(\QQ\cap\left(\QQ+n\right)\right) \\
&\geq \frac{1}{\det \Lat} \vol^{2} \QQ,
\end{align}
where the last inequality is Corollary \ref{main cor}.

To see how inequality  \eqref{easier Van der Corput} of part \eqref{second part of Van der Corput} follows from Theorem \ref{main, for Q},  eq. \eqref{Main2}, we use a similar argument.  Here  $\QQ$ is a convex, centrally symmetric body, so
we have
\begin{equation}\label{copy of eq. from Theorem 3}
    \sum_{n \in \Int(\QQ) \cap \Lat } \vol\Big(\QQ\cap\left(\QQ+2n\right)\Big)
    =\frac{1}{\det \Lat}  \frac{\vol^{2} \QQ}{2^{d}}
    + \frac{2^d}{\det \Lat}
    \sum_{\xi \in \Lat^*\setminus\{0\}}\left|\hat{1}_{\frac{1}{2} \QQ}(\xi)\right|^{2}.
\end{equation}
Now
\begin{align}
   \#\left\{
     \Int(\QQ) \cap \Lat
    \right\}    \vol \QQ 
&=
  \sum_{n \in \Int(\QQ) \cap \Lat } \vol\QQ \\
&\geq 
\sum_{n \in \Int(\QQ) \cap \Lat } \vol\left(\QQ\cap\left(\QQ+2n\right)\right) \\
&\geq 
\frac{1}{\det \Lat}  \frac{\vol^{2} \QQ}{2^{d}},
\end{align}
proving  \eqref{easier Van der Corput}. The last inequality followed from \eqref{copy of eq. from Theorem 3}.




\section{A nonconvex polygon that multi-tiles with multiplicity 2}

\medskip
\begin{example} \label{Nonconvex body 2}
\rm{
 Here we give a nonconvex polygon $\QQ$ (see Figure \ref{crown body})  that multi-tiles with multiplicity $k=2$.  The main point here is to show how the proof of this $2$-tiling follows from  Corollary 
 \ref{main cor}.  
 
 A {\bf nontrivial $k$-tiling} means that there does not exist an $m$-tiling for $m< k$. So here we also have to show that 
 $\QQ$ does not tile with $k=1$ for any lattice.  In other words, we have a non-trivial $2$-tiling.  This phenomenon is in sharp contrast with multi-tiling the plane by using \emph{convex polygons}, because in that convex context the smallest non-trivial multiplicity is $k=5$ 
 \cite{Chuanming1}.

\begin{figure}[ht]
    \centering
\includegraphics[width=5.6cm]{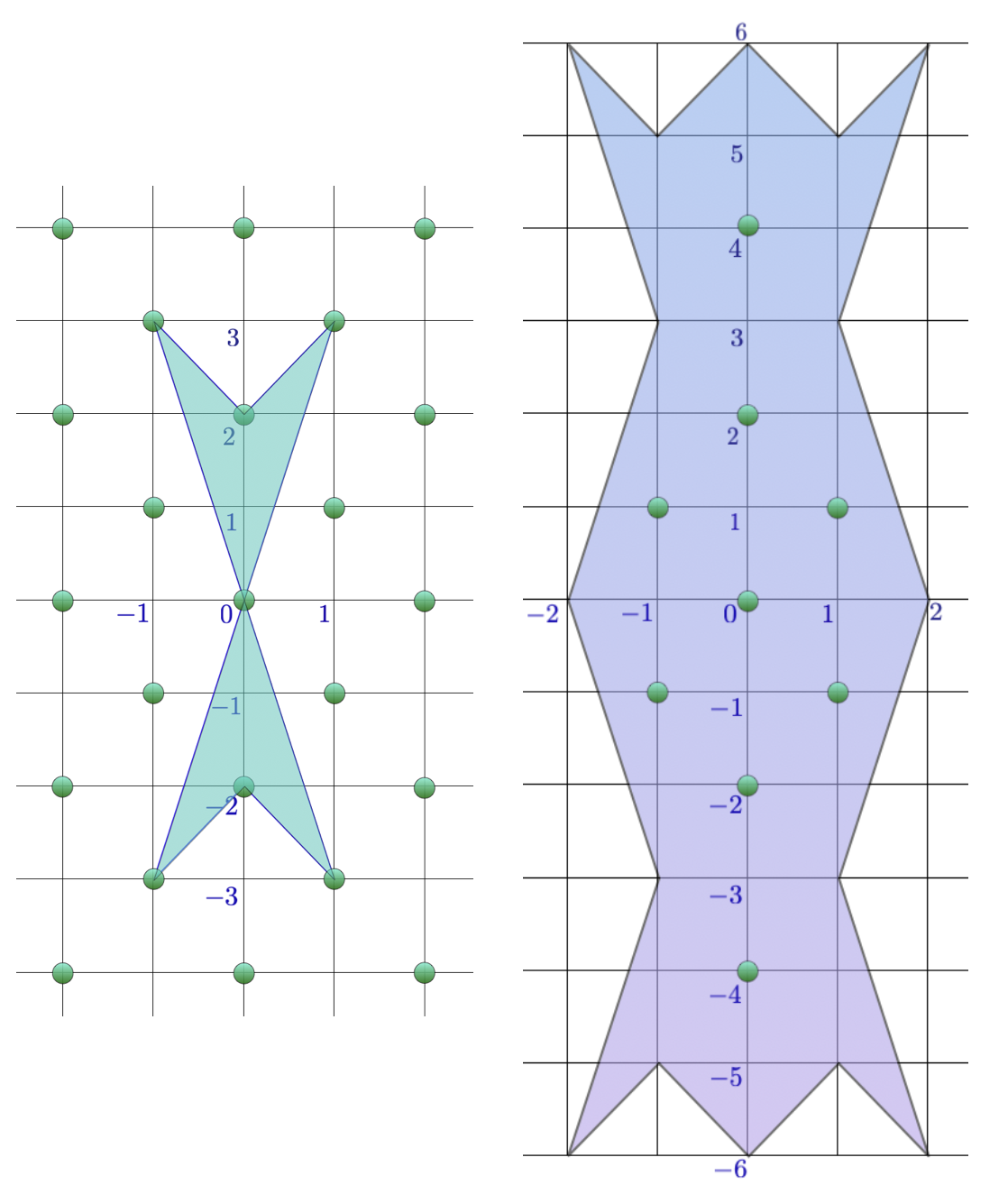}
    \caption{Left: A body $\QQ$, and a lattice $\Lat$ with $\det \Lat = 2$. \newline Right: the difference body $\QQ - \QQ$, with its $9$ interior lattice points of $\Lat$. } \label{crown body}
\end{figure}

 The lattice $\mathcal L$ we chose for this multi-tiling is defined by the basis
 $\left\{ 
 \tvect{1}{1}, \tvect{\ \ 1}{-1} 
 \right\}$, and has index $2$ in $\Z^2$ (also known as the $D_2$ lattice, and drawn with green dots in Figure \ref{crown body}).
\begin{figure}[hbt!]
    \centering
    {{\includegraphics[width=4cm]{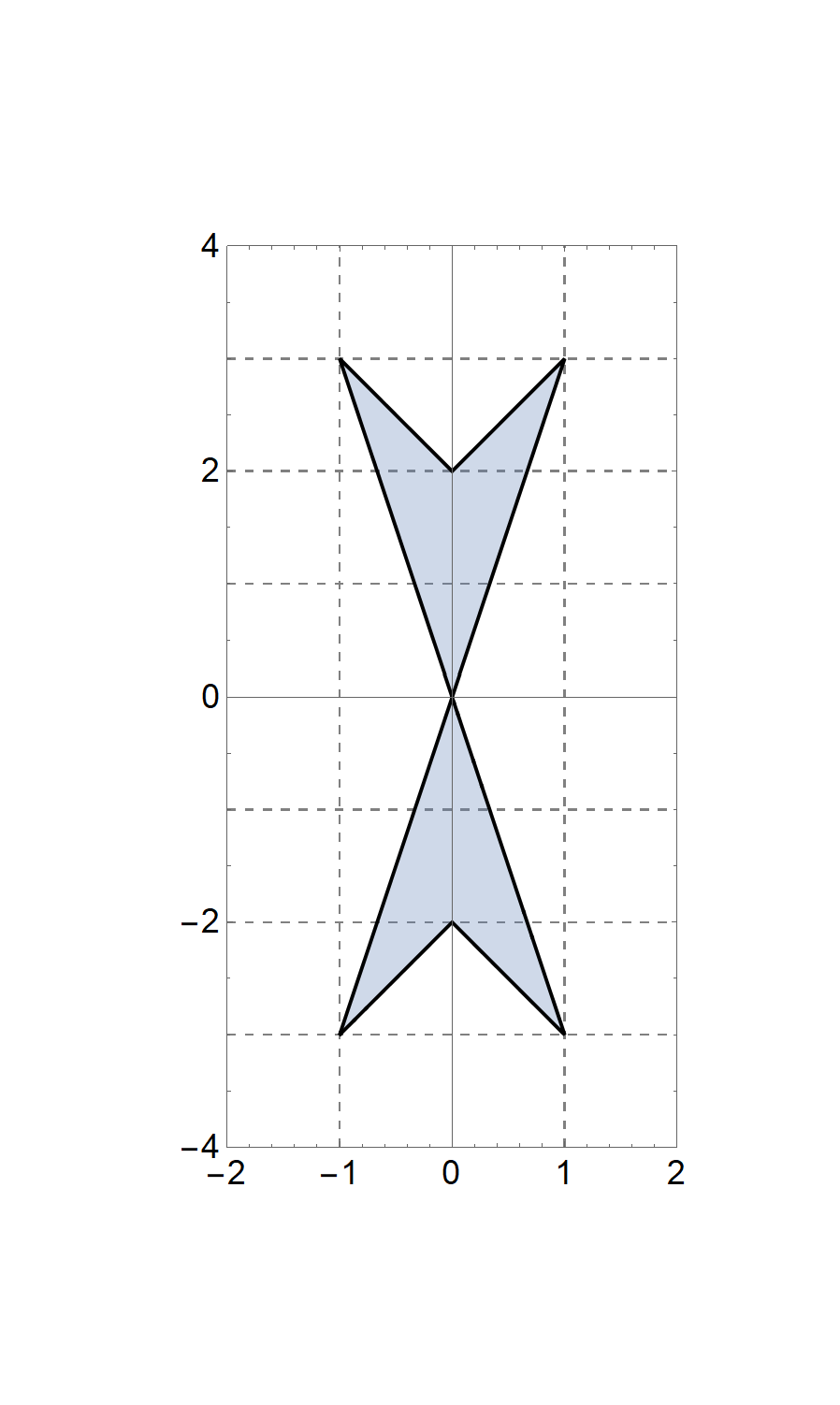} }}
    \qquad
    {{\includegraphics[width=5.5cm]{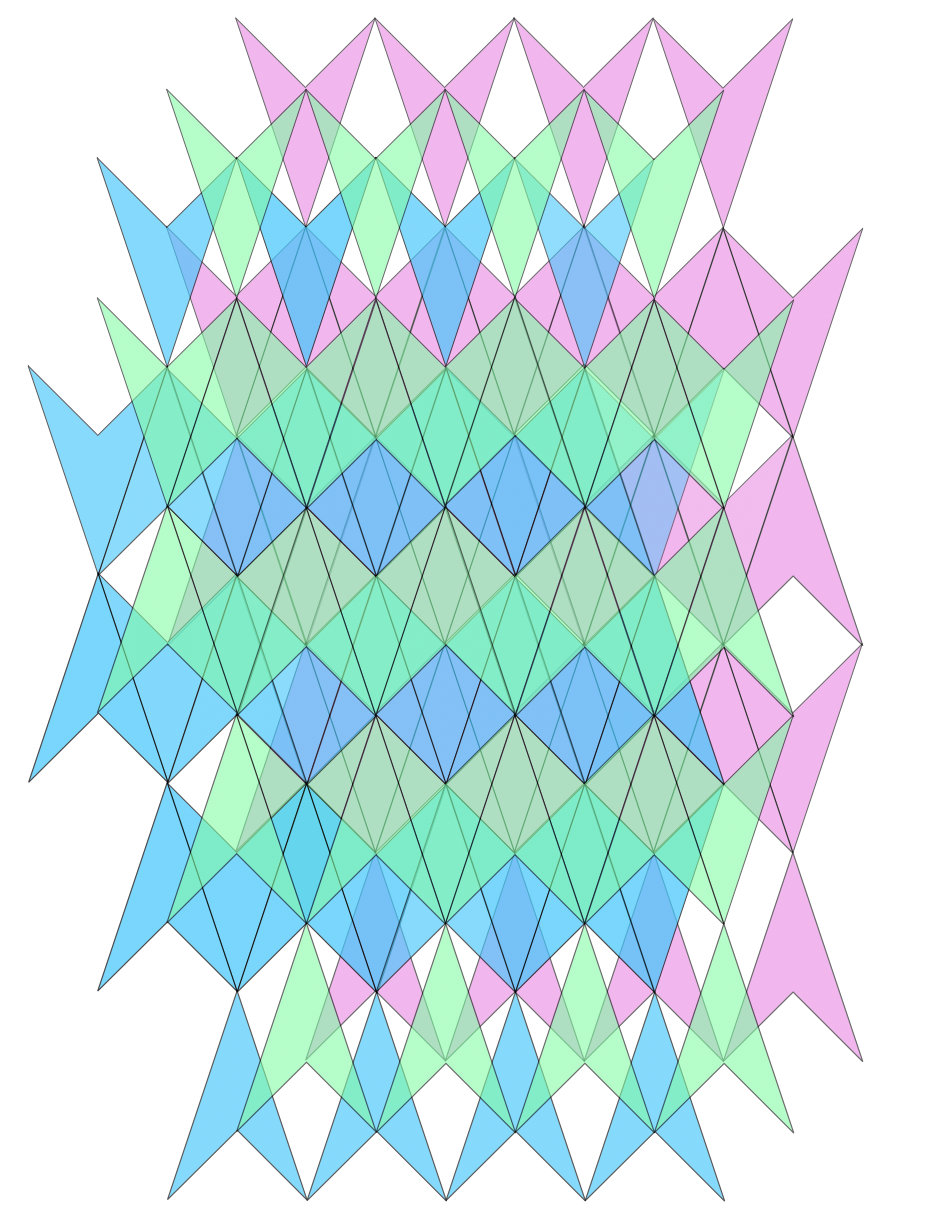} }}
    \caption{$\QQ$ multi-tiles nontrivially, with multiplicity $k=2$.}
    \label{Second nonconvex example}
\end{figure}
By Corollary  \ref{main cor}, to prove that the body $\QQ$ gives a $2$-tiling of the plane, it suffices to show the following equality:
\begin{equation} 
\label{example identity to confirm}
    \sum_{n\in\Lat\cap \interior(\QQ-\QQ)} \vol\left(Q\cap\left(Q+n\right)\right) 
    = 2 \vol Q=8.
\end{equation}
A brute-force computation (using  \rm{Mathematica} \cite{Wolfram}) reveals that the left-hand side of 
\eqref{example identity to confirm} indeed equals $8$, giving us a  computational method of verifying that $\QQ$ indeed $2$-tiles the plane with translations by the lattice $\Lat$.

It is instructive to also see this multi-tiling geometrically, and in
 Figure \ref{Second nonconvex example} we give such a geometric confirmation. 
We first translate the blue collection of translates of $\QQ$ by the vector $\tvect{1}{1}$, to achieve the green collection.  We then translate the same blue collection once again, by the vector $\tvect{2}{2}$, to achieve the pink collection.  Each of the little squares in the middle of the picture show that their color is the overlap of exactly two translates of the blue collection, and hence every point in their interior gets covered exactly twice.  
We therefore obtain a $2$-tiling.

Finally, to see that $\QQ$ does not $1$-tile with any lattice, we suppose to the contrary that it does.   Consulting Figure 
\ref{No 1-tiling}, the translated parallel edge $E_2$, belonging to the purple translate of $\QQ$,  must meet $E_1$, because $E_2$ is the only parallel edge to $E_1$.  Similarly the unique  edge parallel to $E_3$ must be translated to meet $E_3$, as in the right-hand side of Figure \ref{No 1-tiling}.  We arrive at a contradiction, because the two translated copies of $\QQ$, drawn in green and purple, must overlap. 

\begin{figure}[ht]
    \centering
\includegraphics[width=4cm]{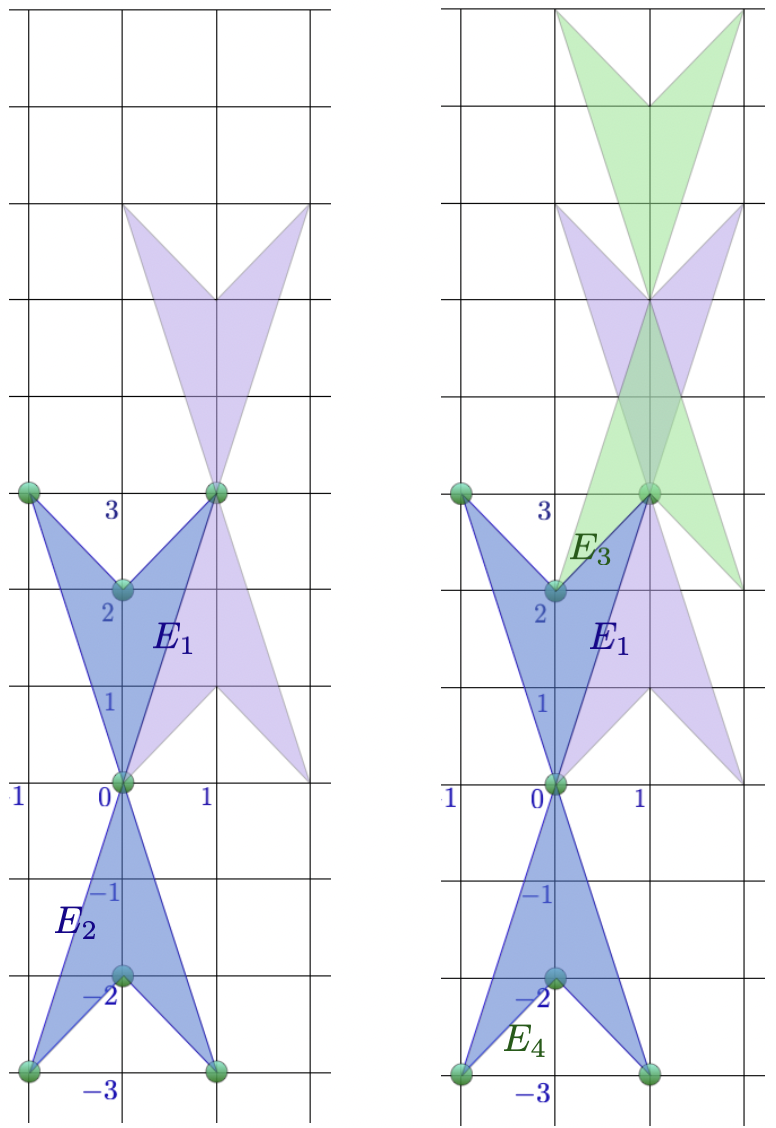}
    \caption{Left: The body $\QQ$ is translated so that edge $E_1$ meets its translated parallel edge $E_2$.  Right: $\QQ$ is translated again, where edge $E_3$ meets its translated parallel edge $E_4$.   } 
    \label{No 1-tiling}
\end{figure}
}
\end{example}


\section{A nonconvex polygon that does not multi-tile}

\medskip
\begin{example} \label{Nonconvex body 1}
\rm{
 Here we illustrate Corollary \ref{main cor} by giving a
 nonconvex polygon does not multi-tile with the integer lattice. Consider the non-convex body $\QQ$ drawn on the left-hand side of  Figure 
\ref{Difference body, take 2} below. 
The difference body of $\QQ$, namely 
$\QQ-\QQ=\{p-q : p,q\in\QQ\}$, is shown on the right-hand side of Figure \ref{Difference body, take 2}, with the red integer points on its boundary.

\begin{figure}[ht]
    \centering
    \includegraphics[width=9cm]{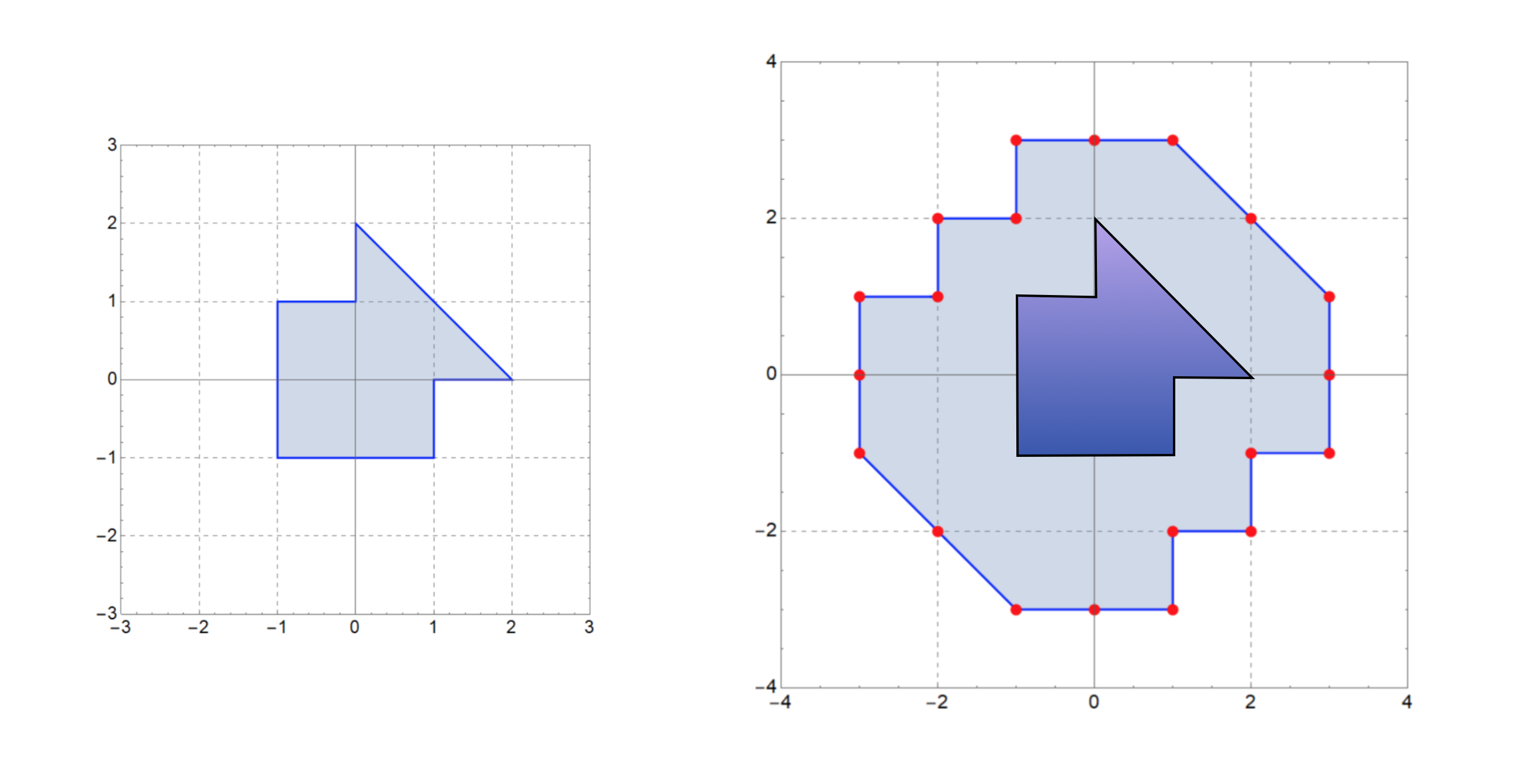}
    \caption{Left:  the nonconvex body $\QQ$ of Example \ref{Nonconvex body 1}.  \newline 
    Right: its difference body $\QQ - \QQ$, with the (red) integer lattice points on its boundary.}
    \label{Difference body, take 2}
\end{figure}

\begin{figure}[ht]
    \centering
    \includegraphics[width=5.5cm]{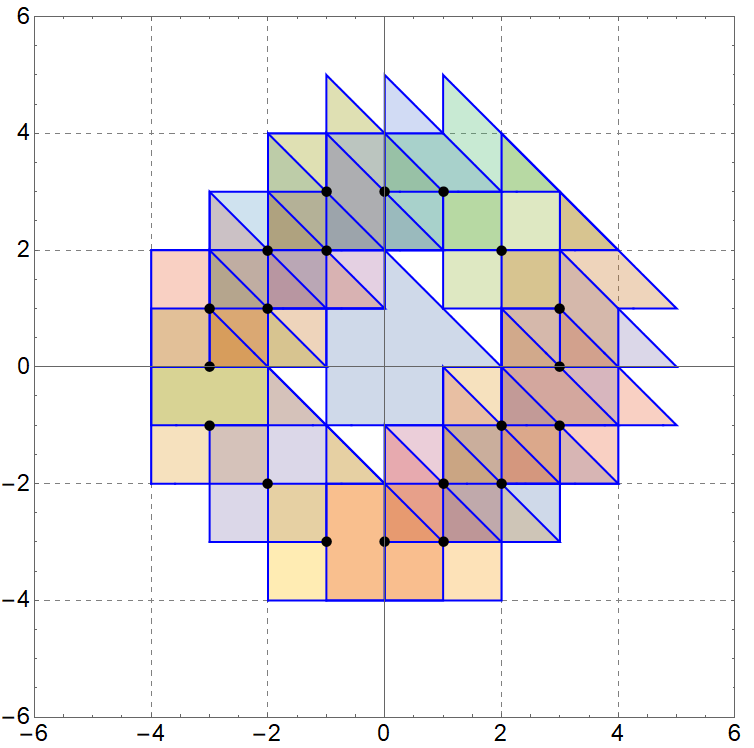}
    \caption{The translated bodies $\QQ+n$ around $\QQ$. As the figure illustrates, the intersections $\QQ\cap(\QQ+n)$, as $n$ varies over the interior points of $\QQ-\QQ$,  only contain boundary points of $\QQ$.}
    \label{Zero volume}
\end{figure}

We will check that here $\QQ$ does not multi-tile with the integer lattice $\Z^2$, by using Corollary \ref{main cor}.  
Namely, we'd like to show that
\begin{equation} \label{Cor 2 for example}
    \sum_{n\in\Lat\cap \interior(\QQ-\QQ)} \vol\left(Q\cap\left(Q+n\right)\right) 
   >  \frac{\vol^2 Q}{\det\Lat}.
\end{equation}
   For this example, $\vol \QQ = 5$, and it turns out the the left-hand-side of
\eqref{Cor 2 for example} equals $26$ (by a brute force computation that used Mathematica), while the right-hand-side equals $ \frac{\vol^2 Q}{\det\Lat} = 25$.  So we've confirmed that indeed we have $26 > 25$ 
in our inequality \eqref{Cor 2 for example}, implying by Corollary \ref{main cor} that  $\QQ$ does not multi-tile with the integer lattice.
}
\hfill $\square$
\end{example}

\section{Proof of 
Theorem \ref{general_case}}
\begin{proof}
\hypertarget{Proof of general_case}
{(of Theorem \ref{general_case})}
We will apply our Poisson summation formula, namely Theorem 
\ref{OurPoissonSummation}, to the function $f*g$.   We first note that because both $f, g, \in L^2(\R^d)$, their convolution $f*g$ is continuous (for example, \cite{Robins}, Exercise $4.19$).  Furthermore, because both $f$ and $g$ are compactly supported, so is $f*g$.  Finally, to satisfy the hypotheses of Theorem 
\ref{OurPoissonSummation} we need to check that $f*g \in L^1(\R^d)$ and that $\widehat{f*g} \in L^1(\R^d)$.  The former condition follows because the convolution of two $L^1$ functions is also in $L^1$ (using Cauchy-Schwartz), and the latter follows because the transform of a compactly supported function (in this case $f*g$) is always in $L^1$.
So by Theorem \ref{OurPoissonSummation}, we have:
\begin{align} 
      \sum_{n\in\Lat}(f*g)(n+x)
&=\frac{1}{\det{\Lat}}
\sum_{m\in\Lat^*}
\hat f(m)
\hat g(m)
e^{2\pi i \langle m, x\rangle},
\end{align}
proving parts \eqref{first part of main result} and \eqref{second part of main result}. 
To prove part \eqref{third part of main result}, we first note that since
 $f,g \in L^2(\R^d)$, we also have
$\hat{f}, \hat{g} \in L^2(\R^d)$.   The Cauchy-Schwarz inequality then gives us:
\begin{equation}
    \int_{\mathbb{R}^{d}}|\reallywidehat{1_A*1_{-B}}(\xi)|\,d \xi= \int_{\mathbb{R}^{d}}|\hat{1}_A(\xi)||\hat{1}_{-B}(\xi)|\,d \xi\leq\left(\int_{\mathbb{R}^{d}}|\hat{1}_A(\xi)|^2\,d \xi\right)^{1/2}\left(\int_{\mathbb{R}^{d}}|\hat{1}_{-B}(\xi)|^2\,d \xi\right)^{1/2}<\infty,
\end{equation}
so that  $\reallywidehat{1_A*1_{-B}}\in L^1(\R^d)$.
We first compute the left hand side of equation \eqref{General Poisson}:

\begin{align*}
    \sum_{n \in \Lat} (1_A*1_{-B})(n+x) &=\sum_{n \in \Lat} \int_{\mathbb{R}^{d}} 1_A(y) 1_{-B}(n+x-y) dy \nonumber\\
    &=\sum_{n \in \Lat} \int_{\mathbb{R}^{d}} 1_{A}(y) 1_{B}(y-n-x) d y\\
    &=\sum_{n \in \Lat} \int_{\mathbb{R}^{d}} 1_{A}(y) 1_{B+n+x}(y) d y\\
    &=\sum_{n \in \Lat} \int_{\mathbb{R}^{d}} 1_{A\cap\left(B+n+x\right)}(y)  d y\\
    &=\sum_{n\in\Lat} \vol\left(A\cap\left(B+n+x\right)\right).
\end{align*}
On the other hand, the right-hand-side of Poisson summation, namely eq. \eqref{General Poisson}, gives us:
\begin{align*}
   \frac{1}{\det \Lat} \sum_{\xi \in \Lat^*} \reallywidehat{1_A*1_{-B}}(\xi)\,e^{2\pi i \langle \xi, x\rangle}
    &=\frac{1}{\det \Lat} \sum_{\xi \in \Lat^*} 
    \hat{1}_{A}(\xi) \hat{1}_{-B}(\xi)\,e^{2\pi i \langle \xi, x\rangle}  
\end{align*}

Summarizing, Theorem \ref{OurPoissonSummation} (our variation of Poisson summation) gives us the required identity
\begin{align*}
\label{Identity1}
    \sum_{n\in\Lat} \vol\left(A\cap\left(B+n+x\right)\right)
    &=
    \frac{1}{\det \Lat} \sum_{\xi \in \Lat^*} 
    \hat{1}_{A}(\xi) \hat{1}_{-B}(\xi)\,e^{2\pi i \langle \xi, x\rangle} \\
&=\frac{1}{\det \Lat} 
    \sum_{\xi \in \Lat^*} 
    \hat{1}_{A}(\xi)
    \overline{\hat{1}_{B}(\xi)}
    \, e^{2\pi i \langle \xi,x \rangle},
\end{align*}
for all $x\in\R^d$.
\end{proof}

\begin{rem}
We note that one direction of the known Theorem \ref{Kol} of Kolountzakis, follows as a Corollary of Bombieri's Theorem.
Namely, suppose that $\QQ$  multi-tiles $\R^d$.  
Then we know by Lemma \ref{simple lemma} that the left hand side of \eqref{Bombieri eq} is constant for all $x$ and therefore by uniqueness of Fourier series, all of the     Fourier coefficients on the right-hand side of \eqref{Bombieri eq} must vanish, except for $\xi=0$.
This guarantees that
Theorem \ref{Kol}, part \ref{Kol, part 2} $\implies$ 
Theorem \ref{Kol}, part \ref{Kol, part 1}, 
and also that $k=\frac{\vol{\QQ}}{\det\Lat}$.
\end{rem}


\section{Proofs of Lemma \ref{first lemma}
and Theorem \ref{main, for Q}}

\begin{proof}(of Lemma \ref{first lemma})
\hypertarget{Proof of first lemma}{We}
observe that the right-hand side of 
\eqref{extending eq. of Bombieri} does not depend on boundary points of $A$ or $B$, which implies that 
\begin{equation}\label{Identity2}
    \sum_{n\in\Lat} \vol\left(A\cap\left(B+n\right)\right)=
    \sum_{n\in\Lat} \vol\left(\Int{A}\cap\left(\Int{B}+n\right)\right).
\end{equation}
But 
\begin{align}
    A\cap\left(B+n\right)&\supset
    \Int{A}\cap\left(\Int{B}+n\right)
    \nonumber\\
    \implies  \vol\left(A\cap\left(B+n\right)\right)
    & \geq\vol\left(\Int{A}\cap\left(\Int{B}+n\right)\right).  
\end{align}
Since the left hand side of equation \eqref{Identity2} is finite, and all terms of both sides of \eqref{Identity2} are non-negative, it follows that
\begin{equation}\label{Identity3}  \vol\left(A\cap\left(B+n\right)\right)=\vol\left(\Int{A}\cap\left(\Int{B}+n\right)\right),
\end{equation}
for each $n\in\Lat$.
By Lemma \ref{basiclemma}, we have
\begin{align}
&\sum_{n\in\Lat} \vol\left(A\cap\left(B+n\right)\right) 
=\sum_{n\in\Lat\cap(A-B)} \vol\left(A\cap\left(B+n\right)\right)\label{Ineq3}\\
    &=\sum_{n\in\Lat\cap\Int{(A-B)}} \vol\left(A\cap\left(B+n\right)\right)+\sum_{n\in\Lat\cap\partial(A-B)} \vol\left(A\cap\left(B+n\right)\right).\nonumber
\end{align}
where in the last step we have an equality since $A-B$ is compact (by assumption), hence closed.
On the other hand, again by Lemma \ref{basiclemma}, we have
\begin{align}
    \sum_{n\in\Lat} \vol\left(\Int{A}\cap\left(\Int{B}+n\right)\right)&=\sum_{n\in\Lat\cap(\Int{A}-\Int{B})} \vol\left(\Int{A}\cap\left(\Int{B}+n\right)\right)\label{Ineq4}\\
    &\leq \sum_{n\in\Lat\cap\Int{(A-B)}} \vol\left(\Int{A}\cap\left(\Int{B}+n\right)\right),\nonumber
\end{align}
where the last inequality is justified by
observing that $\Int{A}-\Int{B}\subset A-B$; therefore if  $\Int{A}-\Int{B}$ is open, then
 $\Int{A}-\Int{B}\subset \Int(A-B)$. 
Combining inequalities (\ref{Ineq3}) and (\ref{Ineq4}), and then using identity (\ref{Identity2}), we get
\begin{equation}\label{Ineq1}
    \sum_{n\in\Lat\cap\Int{(A-B)}} \vol\left(\Int{A}\cap\left(\Int{B}+n\right)\right)\geq\sum_{n\in\Lat\cap\Int{(A-B)}} \vol\left(A\cap\left(B+n\right)\right)+\sum_{n\in\Lat\cap\partial(A-B)} \vol\left(A\cap\left(B+n\right)\right)
\end{equation}
Inserting (\ref{Identity3}) into  (\ref{Ineq1}), we obtain
\begin{equation}
    \sum_{n\in\Lat\cap\partial(A-B)} \vol\left(A\cap\left(B+n\right)\right)\leq 0.
\end{equation}
Since each volume is non-negative, we must have
\begin{equation}
     \vol\left(A\cap\left(B+n\right)\right)=0,
\end{equation}
for each $n\in\Lat\cap\partial(A-B)$ as claimed.
  Part (b)  follows as an immediate consequence, using equation \eqref{Ineq3}.
\end{proof}

\begin{proof}(of Theorem \ref{main, for Q})
\hypertarget{Proof of main, for Q}
To prove part \eqref{part a of Theorem 3}, we will combine Lemma \ref{first lemma}  part \eqref{part b of first lemma}, with Theorem \ref{general_case} part \eqref{third part of main result}, using $x=0$: 
\begin{align*}
\sum_{n\in\Lat\cap \interior(A-B)} \vol\Big(A\cap\left(B+n\right)\Big)&=\sum_{n\in\Lat}
\vol\Big(
A\cap\left(B+n\right)
\Big)\\
&=\frac{1}{\det \Lat} \sum_{\xi \in \Lat^*}
   \hat{1}_{A}(\xi)
 \overline{ \hat{1}_{B}(\xi)}\\
&=\frac{1}{\det \Lat} 
    \vol A \vol B
    + \frac{1}{\det \Lat}
    \sum_{\xi \in \Lat^*\setminus\{0\}}
    \hat{1}_{A}(\xi)
    \overline{\hat{1}_{B}(\xi)}.
\end{align*}

To prove part \eqref{part b of Theorem 3}, we let $A=B=\frac{1}{2}\PP$ in 
\eqref{part a of Theorem 3} above. 
Because $\PP$ is now a centrally symmetric and convex body, we have $A-B=\frac{1}{2}\PP -\frac{1}{2}\PP=\PP$.  Together with $\vol{A}=\vol{B}=\frac{\vol{\PP}}{2^d}$, we obtain: 
\begin{equation}
    \sum_{n \in \Int(\PP) \cap \Lat } \vol\left(\frac{1}{2}\PP\cap\left(\frac{1}{2}\PP+n\right)\right)
    =\frac{1}{\det \Lat}  \frac{\vol^{2} \PP}{2^{2d}}
    + \frac{1}{\det \Lat}
    \sum_{\xi \in \Lat^*\setminus\{0\}}\left|\hat{1}_{\frac{1}{2} \PP}(\xi)\right|^{2}.
\end{equation}
Thus, 
\begin{equation}
    \sum_{n \in \Int(\PP) \cap \Lat } \vol\Big(\PP\cap\left(\PP+2n\right)\Big)
    =\frac{1}{\det \Lat}  \frac{\vol^{2} \PP}{2^{d}}
    + \frac{2^d}{\det \Lat}
    \sum_{\xi \in \Lat^*\setminus\{0\}}\left|\hat{1}_{\frac{1}{2} \PP}(\xi)\right|^{2}.
\end{equation}
\end{proof}

\section{Proofs of Corollaries \ref{main cor}, \ref{multi-tiling Cor}, and 
\ref{max of FT is at the origin}}

The following proof of Corollary \ref{main cor} essentially follows by 
combining Bombieri's Theorem \ref{Bombieri's Theorem}, Kolountakis' Theorem \ref{Kol}, and Lemma \ref{first lemma}.

\begin{proof} 
\hypertarget{proof of main cor}
(of Corollary \ref{main cor})
Since $\QQ \subset \R^d$ is a compact set, equation \eqref{Bombieri eq} gives (with $x=0$): 
\begin{equation} 
\label{inequality for Cor. 2}
   \sum_{n\in\Lat} \vol\left(\QQ\cap\left(\QQ+n\right)\right)=
    \frac{1}{\det \Lat} \sum_{\xi \in \Lat^*}
    \left|\hat{1}_{\QQ}(\xi)\right|^{2}
    \geq
    \frac{1}{\det \Lat}     \left|\hat{1}_{\QQ}(0)\right|^{2} 
    =
\frac{\vol^2 Q}{\det\Lat}.    
\end{equation}
We may rewrite the left-hand side of
\eqref{inequality for Cor. 2}, using 
Lemma \ref{first lemma}, part (b), as follows:
\begin{equation} \label{first part of cor 2}
    \sum_{n\in\Lat\cap \interior(\QQ-\QQ)} \vol\left(Q\cap\left(Q+n\right)\right) 
    \geq  \frac{\vol^2 Q}{\det\Lat}.
\end{equation}
Now we'll show that equality holds in 
\eqref{first part of cor 2} if and only if $\QQ$ multi-tiles, with the lattice $\Lat$. 
First, let's assume that $\QQ$ multi-tiles $\R^d$; that is, Theorem \ref{Kol}, part \ref{Kol, part 2} holds.
Since Theorem \ref{Kol}, part \ref{Kol, part 2}
$\implies$
Theorem \ref{Kol}, part \ref{Kol, part 1},
we have
\begin{equation} 
   \sum_{n\in\Lat\cap \interior(\QQ-\QQ)} \vol\left(\QQ\cap\left(\QQ+n\right)\right)=
    \frac{1}{\det \Lat}\left|\hat{1}_{\QQ}(0)\right|^{2}=\frac{\vol^2(\QQ)}{\det \Lat}.
\end{equation}
Conversely, if equality holds in 
\eqref{first part of cor 2}, then
\begin{equation} 
   \sum_{n\in\Lat\cap \interior(\QQ-\QQ)} \vol\left(\QQ\cap\left(\QQ+n\right)\right)=\frac{\vol^2(\QQ)}{\det \Lat},
\end{equation}
and by 
\eqref{Bombieri eq} with $x=0$ we therefore must have $\left|\hat{1}_{\QQ}(\xi)\right|^{2}=0$ for all $\xi\in\Lat^*$, excluding the origin. This means that 
Theorem \ref{Kol}, part \ref{Kol, part 1} holds, and its equivalence with 
Theorem \ref{Kol}, part \ref{Kol, part 2} tells us that
$\QQ$ multi-tiles $\R^d$ with the lattice $\Lat$.
\end{proof}

\bigskip
\begin{proof}
\hypertarget{proof of multi-tiling Cor}
(of Corollary \ref{multi-tiling Cor})
(a) $\Rightarrow$ (b). Taking $A=B$ in  Theorem \ref{main, for Q} part (a) we have
\begin{align*}
   \sum_{n \in \Int(A-A) \cap \Lat } \vol\Big(A\cap\left(A+n\right)\Big)
    &=\frac{1}{\det \Lat} 
    \vol A \vol B
    + \frac{1}{\det \Lat}
    \sum_{\xi \in \Lat^*\setminus\{0\}}
    \hat{1}_{A}(\xi)
    \overline{\hat{1}_{B}(\xi)}\\
    &=\frac{1}{\det \Lat} 
    \vol^2 A 
    + \frac{1}{\det \Lat}
    \sum_{\xi \in \Lat^*\setminus\{0\}}
    |\hat{1}_{A}(\xi)|^2,
\end{align*}
implying that 
\begin{equation}
    \frac{1}{\det \Lat}
    \sum_{\xi \in \Lat^*\setminus\{0\}}
    |\hat{1}_{A}(\xi)|^2=0.
\end{equation}
Hence $\hat{1}_{A}(\xi)=0$ for all $\xi \in \Lat^*\setminus\{0\}$. By Theorem \ref{Kol} we conclude that  $A$ $k$-tiles $\R^d$ by translations with $\Lat$.

To prove (b) $\Rightarrow$ (c), we use
Theorem \ref{Kol}.  Because $A$ $k$-tiles $\R^d$, $\hat{1}_{A}(\xi)=0$ for all $\xi \in \Lat^*\setminus\{0\}$. Again, by Theorem \ref{main, for Q} part (a), we get 
\begin{align*}
   \sum_{n \in \Int(A-A) \cap \Lat } \vol\Big(A\cap\left(A+n\right)\Big)
    &=\frac{1}{\det \Lat} 
    \vol A \vol B
    + \frac{1}{\det \Lat}
    \sum_{\xi \in \Lat^*\setminus\{0\}}
    \hat{1}_{A}(\xi)
    \overline{\hat{1}_{B}(\xi)}\\
    &=\frac{1}{\det \Lat} 
    \vol A\vol B 
    + \frac{1}{\det \Lat}
    \sum_{\xi \in \Lat^*\setminus\{0\}}
    0\cdot\overline{\hat{1}_{B}(\xi)}\\
    &=\frac{1}{\det \Lat} 
    \vol A\vol B ,
\end{align*}
and we are done.

Finally, to prove (c) $\Rightarrow$ (a), we just take $A=B$ in part (c) 
(eq. \eqref{part (c) of Cor 2}). 
\end{proof}


\bigskip
\begin{proof}
\hypertarget{proof of max of FT is at the origin}
(of Corollary \ref{max of FT is at the origin})
This proof  follows from \eqref{Bombieri eq}  by observing that the hypotheses guarantee that the left-hand side of equation \eqref{alternate main Theorem} vanishes. Namely, we have:
\begin{equation}
  0= \sum_{n\in\Lat} \vol\left(\QQ\cap\left(\QQ+n+ x\right)\right)=
    \frac{1}{\det \Lat}  
    \sum_{\xi \in \Lat^*}
    \left|\hat{1}_{\QQ}(\xi)\right|^{2}
    \cos{(2\pi\langle \xi, x \rangle)}.
\end{equation}
Since $\hat{1}_{\QQ}(0) = \vol \QQ$, we're done.
\end{proof}
\section{Proof of Theorem 
\ref{application: arithmetic combintorics 1} }

\begin{proof}
\hypertarget{first theorem for integer sets}
To prove part \eqref{first eq. of AC proof}, we simply substitute $\vol\left(A\cap\left(A+n\right)\right) = |F\cap (F+n)| \varepsilon^d$ and $\vol A = \varepsilon^d |F|$
into equation \eqref{special discrete case 1}, and divide both sides of the ensuing equation by $\varepsilon^d$.

To prove part \eqref{second eq. of AC proof}, we merely
have to simplify the expression
$\left | \hat{1}_{A}(\xi)\right |^2$. 
The main point here is that we have the disjoint union
\[
A  = \bigcup_{a\in F} \left( \interior\square + a \right).
\]
In other words, because $0<\epsilon \leq 1$,
the interiors of the 
translated $\varepsilon$-cubes are pairwise disjoint.
Using this disjoint union, we compute:
\begin{align}
\hat{1}_{A}(\xi) &= \sum_{a\in F}
\hat 1_{\square+a}(\xi) 
=\sum_{a\in F} 
\hat 1_{\square}(\xi)
   e^{2\pi i\langle  \xi, a\rangle} 
=   
\hat 1_{\square}(\xi)
   \sum_{a\in F}  e^{2\pi i\langle  \xi, a\rangle} \\
&=   \varepsilon^d 
\prod_{k=1}^d
  \sinc(\pi \varepsilon \xi_k)
   \sum_{a\in F} e^{2\pi i\langle  \xi, a\rangle},
\end{align}
where we've  used the standard $\sinc$ formula 
$\hat 1_{\left[ -\tfrac{1}{2}, \tfrac{1}{2}\right]^d}(\xi)=\prod_{k=1}^d \sinc( \pi \xi_k )$  
for the Fourier transform of the unit cube (see for example \cite{Robins}), as well as the general property $\hat 1_{\varepsilon S}(\xi) =
\varepsilon^d \hat 1_{S}(\varepsilon \xi)$. 
\end{proof}

\section{Proof of Theorem 
\ref{tiling the integer lattice with a finite set of integers} }

\begin{proof}
\hypertarget{main theorem for integer lattices}
To prove \eqref{part a of lattice tilings}  $\implies$ 
\eqref{part b of lattice tilings}, we note that
because the set $F$ multi-tiles $\Z^d$,  its 1-thickening $A$ multi-tiles $\R^d$ so we can apply the implication 
$(a) \implies (b)$ of Corollary \eqref{multi-tiling Cor}  for the compact set $A$: 
\begin{equation} 
    \sum_{n\in\Lat\cap \interior({A}-{A})} \vol\Big({A}\cap\left({A}+n\right)\Big) 
   = \frac{\vol^2 A}{\det\Lat}.
\end{equation}
But, $\vol\Big(A\cap\left(A+n\right)\Big) = |F\cap (F+n)| \varepsilon^d=|F\cap (F+n)|$, and $\vol A = \varepsilon^d |F|=|F|$, proving part \eqref{part b of lattice tilings}. Similarly, the implication $(b) \implies (a)$ from Corollary \eqref{multi-tiling Cor} proves that
\eqref{part b of lattice tilings}  $\implies$ 
\eqref{part a of lattice tilings} here.

To prove that \eqref{part b of lattice tilings}
$\implies$ \eqref{part c of lattice tilings}, 
we recall Corollary \eqref{application: arithmetic combintorics 1}, part \eqref{second eq. of AC proof}, with $\varepsilon=1$:
\begin{equation} \label{proof of b implies c} 
      \sum_{n\in (F-F)\cap \Lat} 
   |F\cap (F+n)| 
  =\frac{1}{\det \Lat} |F|^2
  +
 \frac{1}{ \det \Lat} 
  \sum_{\xi \in \Lat^*\setminus 0} 
  \left(
\prod_{k=1}^d \sinc^2(\pi  \xi_k)
 \left|
 \sum_{n \in F}  e^{2\pi i\langle  \xi, n\rangle}
 \right|^2
 \right),
 \end{equation}
 which holds for any finite set $F$ of integer points.
 Therefore, the assumption of 
 \eqref{part b of lattice tilings} is equivalent to the statement: 
 \begin{equation}\label{vanishing dual lattice sum} 
\sum_{\xi \in \Lat^*\setminus 0} \,
\prod_{k=1}^d \sinc^2(\pi  \xi_k)
 \left|
 \sum_{n \in F}  e^{2\pi i\langle  \xi, n\rangle}
 \right|^2=0.
 \end{equation}
Because all of the summands in the latter sum are nonnegative, equation \eqref{vanishing dual lattice sum} is equivalent to
\begin{equation}
    \prod_{k=1}^d \sinc^2(\pi  \xi_k)
 \left|
 \sum_{n \in F}  e^{2\pi i\langle  \xi, n\rangle}
 \right|^2=0,
\end{equation}
for each $\xi \in \Lat^*\setminus 0$.
We recall  that 
$\frac{\sin(\pi \xi_k)}{\pi \xi_k}=0 \iff
\xi_k \in \Z$; so we see that for each nonzero $\xi \in \Lat^*$ either one of its coordinates is an integer, or else 
$\sum_{n \in F}  e^{2\pi i\langle  \xi, n\rangle}=0$.
This last observation is precisely the assertion of 
part  \eqref{part c of lattice tilings}, and hence part \eqref{part c of lattice tilings} is equivalent to 
part  \eqref{part b of lattice tilings}.
\end{proof}

\section{More examples}

\begin{example}
\rm{
While it is tempting to assert that $\interior(Q-Q) = \interior Q - \interior Q$, for any compact set $Q$, we give a rather extreme counterexample to this claim. 
Consider the usual cantor set $Q\subset [0,1]$, whose interior is known to be the empty set.  It is known \cite{Kraft} that its difference body satisfies the surprising identity
$Q-Q = [-1, 1]$, implying that $\interior(Q-Q)=(-1, 1)$.  
However, $\interior Q = \varnothing$, so we find that 
$\interior Q - \interior Q = \varnothing \not=
\interior(Q-Q)$.
This example shows why the analysis in Lemma \ref{first lemma} is necessary.
}
\hfill $\square$
\end{example}


\medskip
\begin{example}
\rm{
Suppose we have a convex, centrally symmetric body $Q\subset \R^d$ that $k$-tiles $\R^d$ with a lattice $\Lat$.  
Consequently, we know that $\QQ - \QQ= 2\QQ$.  
Here we use Theorem \ref{main, for Q} to show that if $\QQ$ $k$-tiles, and enjoys the property that
$\QQ - \QQ$
contains exactly $3$ lattice points of $\Lat$ in its interior, say $0, n_0, -n_0$,
then either $k=1$ or $k=2$.

By Theorem \ref{main, for Q} , we have:
\begin{align*}
  \frac{\vol^{2} \QQ}{\det \Lat}
&=
\sum_{n \in \text{int}(\QQ-\QQ) \cap \Lat } \vol\left(\QQ\cap\left(\QQ+n\right)\right) \\
&= \vol \QQ + \vol\left(\QQ\cap\left(\QQ+n_0\right)\right)
+ \vol\left(\QQ\cap\left(\QQ-n_0\right)\right)
\\
&= \vol \QQ + 2\vol\left(\QQ\cap\left(\QQ+n_0\right)\right).
\end{align*}
By Theorem \ref{Kol}, we  also know that
$k = \frac{\vol \QQ}{\det \Lat}$, so that 
\begin{align} \label{first k constraint}
  k \vol \QQ 
&= \vol \QQ + 2\vol\left(\QQ\cap\left(\QQ+n_0\right)\right).
\end{align}
Since $\QQ\cap(\QQ+n_0)\subset\QQ$ we must have
\begin{align*}
  k \vol \QQ 
&= \vol \QQ + 2\vol\left(\QQ\cap\left(\QQ+n_0\right)\right)< 3\vol\QQ,
\end{align*}
giving us 
\[
k \leq 2,
\]
as was claimed.
As an aside, together with equation \eqref{first k constraint}, we now also have
\begin{align}
  k &= 1 + 2\frac{
\vol\left(\QQ\cap\left(\QQ+n_0\right)\right)
}
{
\vol \QQ 
}\leq2.
\end{align}
and therefore
\begin{equation} \label{a little claim}
\vol\left(\QQ\cap\left(\QQ+n_0\right)\right)
 =\frac{1}{2}\vol \QQ \quad\text{or}\quad \vol\left(\QQ\cap\left(\QQ+n_0\right)\right)=0.
\end{equation}
\hfill $\square$
}
\end{example}

\begin{example} \label{zeta(2) example}
\rm{
The classical identity of Euler, namely 
$
 \zeta(2):=1+\frac{1}{2^2}+\frac{1}{3^2}+\frac{1}{4^2}+\cdots=\frac{\pi^2}{6}
$
can be retrieved from 
Bombieri's theorem (eq. \eqref{Bombieri eq}) by using the integer lattice $\Lat=\Z^2$ together with the triangle 
$\QQ:=\triangle$, defined by the vertices 
$v_0:= \tvect{0}{0},   v_1:= \tvect{1}{0}$, and 
 $v_2:= \tvect{0}{1}$,
as shown in Figure \ref{Pi squared over 6}.
\begin{figure}[ht]
    \centering
    \includegraphics[width=5cm]{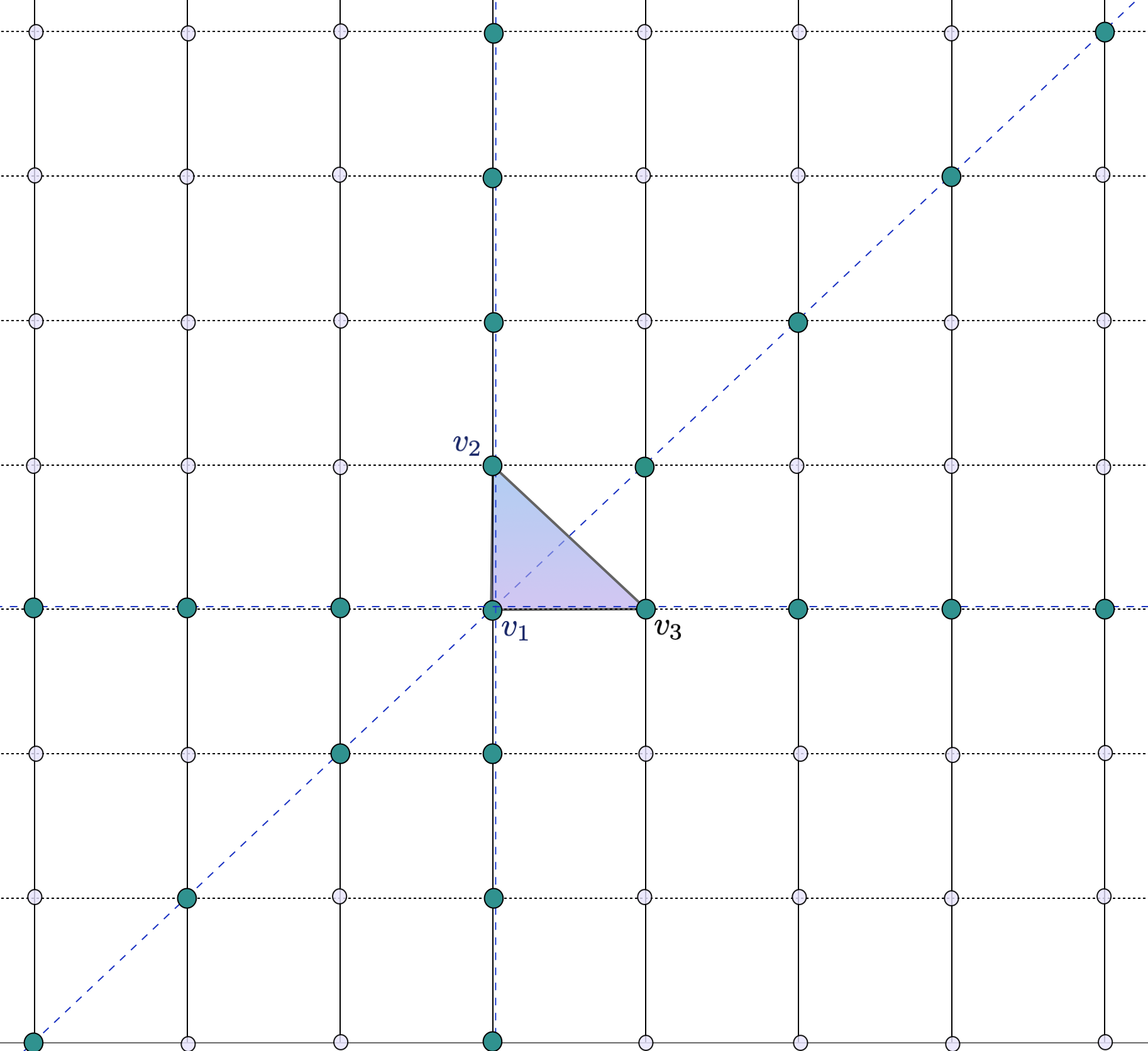}
    \caption{A triangle $\triangle$, and the three discrete lines of integer points (drawn in dark green) where $\hat{1}_{\triangle}(\xi)$ does not vanish.}
    \label{Pi squared over 6}
\end{figure}
On the left-hand side of equation 
\eqref{Bombieri eq}, we have (with $x=0$):
\[
\sum_{n\in\Z^2} \vol\left(\triangle\cap\left(\triangle+n\right)\right)=\vol\triangle=\frac{1}{2}.
\]
On the right-hand side of \eqref{Bombieri eq}
we can compute directly the summands 
$ \left|\hat{1}_{\triangle}(\xi)\right|^{2}$, using the following known formula.  
\begin{lem}\label{unecessary coordinates} For the standard $d$-dimensional simplex 
$\triangle_d:=  \{ x\in \R_{\geq 0}^d \mid  
x_1 + \cdots + x_d \leq 1\} \subset \R^d$,  
its Fourier transform vanishes for all integer points $\xi$ with the property that all of their coordinates are nonzero and distinct from each other. 
\end{lem}
\begin{proof}
The Fourier transform of $\triangle_d$ is given by:
\[
\hat{1}_{\triangle_d}(\xi)
:= \int_{\triangle_d}  
e^{ -2\pi i \langle \xi,  x\rangle} dx
=\frac{1}{(2\pi i)^d}
\sum_{i=0}^{d} \frac{ e^{-2\pi i \xi_{i}}}
{\prod_{0\leq j \leq d,  \, j \neq i }
\left(\xi_{i}-\xi_{j}\right)},
\]
where we let $\xi_i\neq \xi_j$ 
for every $i\neq j$, and where 
we've defined $\xi_0:= 0$
(\cite{Barvinok}, or \cite{Lasserre1}, Lemma 21).  
If we restrict $\xi$ to be an integer point, then for each coordinate 
$\xi_i \in \Z \, (i\geq1)$,  we have  
$e^{-2\pi i \xi_{i}}=1$.  By the Lagrange interpolation polynomial applied to the points $(\xi_i, 1)\in \R^2$, we have $\hat{1}_{\triangle_d}(\xi)=0$.  The reason for this vanishing is that these points lie on a horizontal line - the constant polynomial. 
Consequently, all coefficients  in the interpolating formula vanish except the constant coefficient, which equals 1 by uniqueness of the  polynomial with degree less than or equal to $d+1$.
\end{proof}


From the discussion above, it is sufficient to restrict attention to only those integer points $\xi \in \Z^d$ which have at least one vanishing coordinate, or at least two equal coordinates.
When $d=2$, the only integer vectors $\tvect{\xi_1}{\xi_2}\in\Z^2$ for which $\hat{1}_{\triangle_2}(\xi)$ possibly does not vanish are those that are orthogonal to the sides of $\triangle_2$, namely the integer points belonging to
the family
$\mathcal{F}=\left\{\tvect{0}{k} \mid k \in \Z \right\}
\cup
\left\{\tvect{k}{0} \mid k \in \Z \right\}
\cup
\left\{\tvect{k}{k} \mid k \in \Z \right\}$, as drawn in Figure \ref{Pi squared over 6}. 
By symmetry,  
$\left|\hat{1}_{\triangle_2} \tvect{0}{k} \right|^2
=\left|\hat{1}_{\triangle_2}\tvect{k}{0}\right|^2
$, for all $k\in \Z$.  We compute, for each $k \in \Z-\{0\}$: 
\begin{align*}
    \left|\hat{1}_{\triangle_2}\tvect{0}{k}\right|^2 &=\left|\int_{0}^1\int_{0}^{-x_2+1}  e^{-2 \pi i kx_1} dx_1dx_2\right|^2  \\
    &=\frac{1}{4\pi^2 k^2}\left|\int_{0}^1 (e^{2 \pi i kx_2}-1) dx_2\right|^2\\
    &=\frac{1}{4\pi^2 k^2}.
\end{align*}
A similar computation gives $\left|\hat{1}_{\triangle_2}\tvect{k}{k}\right|^2
=\frac{1}{4\pi^2 k^2}.$ 
By Bombieri's theorem (eq. \eqref{Bombieri eq}), we have:
\begin{align*}
    \frac{1}{2}=\vol\triangle_2
    &=\frac{1}{\det\Z^2}\sum_{\xi\in\mathcal{F}}\left|\hat{1}_{\triangle_2}(\xi)\right|^2\\
    &=\left|\hat{1}_{\triangle_2}(0)
    \right|^2+3\sum_{k\in\Z-\{0\}}\frac{1}{4\pi^2 k^2}\\
    &=(\vol\triangle_2)^2+\frac{3}{4\pi^2}
    \sum_{k\in\Z-\{0\}}\frac{1}{k^2}\\
    &=\frac{1}{4}+\frac{3}{2\pi^2}
    \sum_{k=1}^{\infty}\frac{1}{k^2},
\end{align*}
giving the required classical identity
$\sum_{k=1}^{\infty}\frac{1}{k^2}=\frac{\pi^2}{6}$.
}
\hfill $\square$
\end{example}

\begin{figure}[ht]
    \centering
    \includegraphics[width=6cm]{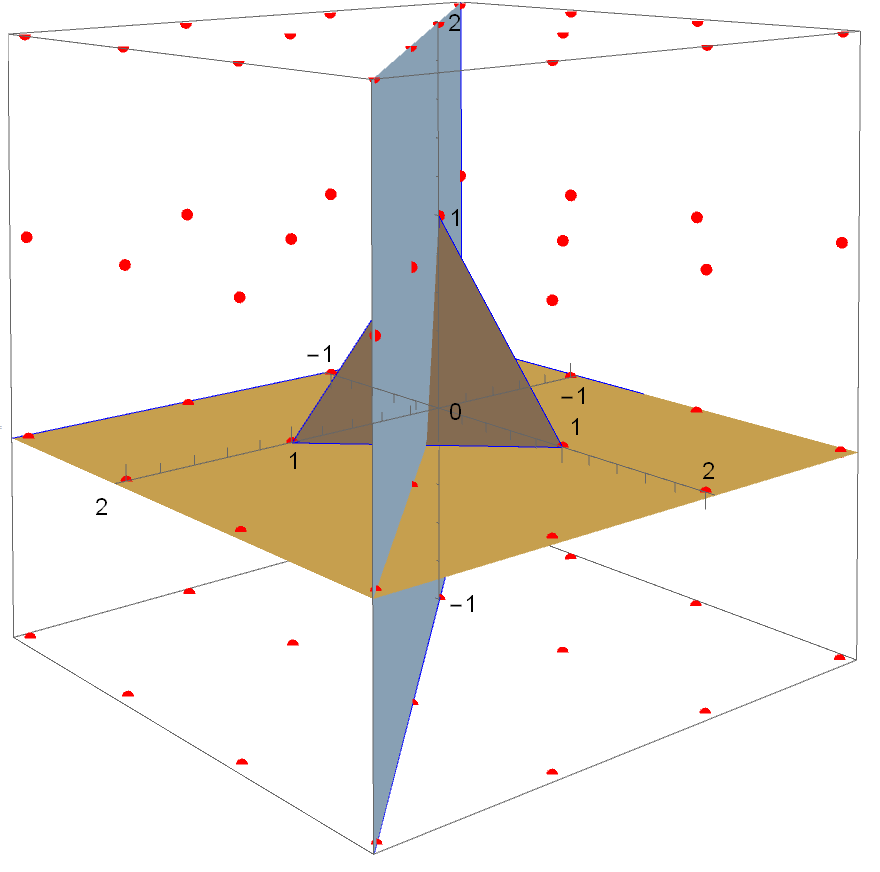}
    \caption{The Tetrahedron $\triangle_3$, the integer lattice $\Z^3$, and the hyperplanes where $\hat{1}_{\triangle_3}(\xi)
    \neq0$. There are four more analogous hyperplanes with this property.}
    \label{Pi fourth power over 90}    
\end{figure}

\begin{example}
\label{odd zeta values via our formula}
\rm{
Here we show that it is possible to recover the values of the Riemann zeta function at all positive even integers, namely $\zeta(2d)$, by using \eqref{Bombieri eq}.  This approach extends Example \ref{zeta(2) example}, and we begin by
considering the integer lattice, together with the standard simplex  $\QQ:=\triangle_d$. 
We will compute
$ \zeta(4):=\sum_{n=1}^\infty \frac{1}{n^4}$.

First, we have
$ \sum_{n\in\Z^3} \vol\left(\triangle_3\cap\left(\triangle_3+n\right)\right)=\vol\triangle_3=\frac{1}{6}$, 
The crucial ingredient on the right-hand side of equation \eqref{Bombieri eq} is the transform 
$
\left|\hat{1}_{\triangle_3}(\xi)\right|^{2}=\left|\int_{\triangle_3}  e^{-2 \pi i\langle\xi, x\rangle} d x\right|^2
$.
The only $\xi=(\xi_1, \xi_2, \xi_3)\in\Z^3$ for which the latter integral does not vanish are those integer points lying on $6$ particular hyperplanes, by  Lemma \ref{unecessary coordinates}.
More precisely, these are exactly the integer points in the family 
\begin{equation*}
    \mathcal{F}=\{(0,k,k),(k,0,k),(k,k,0),(k,k,l), (k,l,k), (l,k,k)\}_{k,l\in\Z}.
\end{equation*}
After some computations we obtain
\begin{align*}
    &\left|\hat{1}_{\triangle_3}(0,k,l)\right|^2=\frac{1}{16\pi^4}\cdot\frac{1}{k^2l^2} \text{  for } k,l\neq0 \text{ and } k\neq l,\\
    &\left|\hat{1}_{\triangle_3}(0,k,k)\right|^2=\frac{1}{16\pi^4}\cdot\frac{4}{k^4} \text{  for } k\neq0,\\
    &\left|\hat{1}_{\triangle_3}(k,k,l)\right|^2=\frac{1}{16\pi^4}\cdot\frac{1}{k^2(k-l)^2} \text{  for } k,l \neq0 \text{ and } k\neq l,\\
    &\left|\hat{1}_{\triangle_3}(0,0,k)\right|^2=\left|\hat{1}_{\triangle_3}(k,k,k)\right|^2=\frac{1}{16\pi^4}\cdot\left(\frac{1}{k^4}+\frac{\pi^2}{k^2}\right)\text{  for } k\neq0.
\end{align*}

By symmetry, all permutations of the coordinates yield the same result for each computation above. By \eqref{Bombieri eq} we have
\begin{align*}
    \frac{1}{6}&=\vol\triangle_3=\frac{1}{\det\Z^3}\sum_{\xi\in\mathcal{F}}\left|\hat{1}_{\triangle_3}(\xi)\right|^2=\left|\hat{1}_{\triangle_3}(0,0,0)\right|^2+3\sum_{k,l\neq0, k\neq l}\left|\hat{1}_{\triangle_3}(0,k,l)\right|^2+3\sum_{k,l\neq0, k\neq l}\left|\hat{1}_{\triangle_3}(k,k,l)\right|^2+\\
    &\hspace{5.5cm}+3\sum_{k\neq0 }\left|\hat{1}_{\triangle_3}(0,0,k)\right|^2+\sum_{k\neq0}\left|\hat{1}_{\triangle_3}(k,k,k)\right|^2+3\sum_{k\neq0}\left|\hat{1}_{\triangle_3}(0,k,k)\right|^2
    \\
    &=(\vol\triangle_3)^2+\frac{1}{16\pi^4}\Bigg(3\sum_{k,l\neq0, k\neq l}\frac{1}{k^2l^2}+3\sum_{k,l\neq0, k\neq l
    }\frac{1}{k^2(k-l)^2}+3\sum_{k\neq0}\left(\frac{1}{k^4}+\frac{\pi^2}{k^2}\right)+\sum_{k\neq0
    }\left(\frac{1}{k^4}+\frac{\pi^2}{k^2}\right)+3\sum_{k\neq0}\frac{4}{k^4}\Bigg)\\
    &=\frac{1}{36}+\frac{1}{16\pi^4}\Bigg(3\sum_{k,l\neq0}\frac{1}{k^2l^2}-3\sum_{k\neq0}\frac{1}{k^4}+3\sum_{k\neq0, k\neq l
    }\frac{1}{k^2(k-l)^2}-3\sum_{k\neq0 }\frac{1}{k^4}+4\sum_{k\neq0}\left(\frac{1}{k^4}+\frac{\pi^2}{k^2}\right)+3\sum_{k\neq0}\frac{4}{k^4}\Bigg)\\
    &=\frac{1}{36}+\frac{1}{16\pi^4}\Bigg(24\sum_{k,l>0}\frac{1}{k^2l^2}+20\sum_{k>0 }\frac{1}{k^4}+8\sum_{k>0 }\frac{\pi^2}{k^2}\Bigg)\\
    &=\frac{1}{36}+\frac{1}{16\pi^4}\Bigg(24\left(\frac{\pi^2}{6}\right)^2+20\sum_{k>0 }\frac{1}{k^4}+8\pi^2\cdot\frac{\pi^2}{6}\Bigg)\\
    &=\frac{1}{36}+\frac{1}{16\pi^4}\Bigg(2\pi^4+20\sum_{k>0 }\frac{1}{k^4}\Bigg)\\
    &=\frac{11}{72}+\frac{5}{4\pi^4}\sum_{k>0 }\frac{1}{k^4},
\end{align*}
giving us the well-known identity  $\sum_{k=1}^{\infty}\frac{1}{k^4}=\frac{\pi^4}{90}$.
}
\hfill $\square$
\end{example}

\bigskip
\section{Further remarks}

Here we mention some open problems, for future research.
\begin{problem}
Extend the notion of Poisson summation friendly sets to include unbounded sets, and classify them. 
\end{problem}
An answer to this problem would enable a generalization of the current results to functions whose support is unbounded.  
Regarding multi-tiling bodies, the following simple-sounding question is still open.
\begin{problem} \label{k-tiling problem}
If $\QQ$ is a convex body that $k$-tiles $\R^d$ (nontrivially), then $k \not=2$. 
\end{problem}
Problem \ref{k-tiling problem} is known to be true in dimensions $2$ and $3$, due to the recent work of \cite{Chuanming1}, but it is open for $d \geq 4$.

\begin{problem} \label{Application: Matheron's conj}
Apply the methods and results herein to study Matheron's Conjecture \ref{conj:Matheron}. 
\end{problem}

\begin{problem} \label{zeta function problem}
Using Theorems \ref{Bombieri's Theorem}  and \ref{main, for Q}, with some appropriate Poisson summation friendly set $\QQ$ (possibly nonconvex),
find a geometric interpretation for  $\zeta(3), \zeta(5)$, etc.
\end{problem}

We note that there are already some fascinating geometric interpretations for the odd special values $\zeta(2n+1)$, following the work of Ed Witten on volumes of certain moduli spaces of flat connections (\cite{Witten}, eq. 4.93).  In Example  \ref{odd zeta values via our formula}, we picked simplices, and we retrieved the even values of $\zeta(s)$ from a spectral representation of their volumes.  Perhaps more can be done for Problem \ref{zeta function problem} by picking more complex sets and computing their volumes in eq.  \eqref{Bombieri eq}, or eq. \eqref{extending eq. of Bombieri}. 


\bigskip \noindent
{\bf Acknowledgement}. \ 
The authors would like to thank Mihalis Kolountzakis for useful remarks.

\newpage
\appendix

\section{Some useful and known lemmas}
\label{Known lemmas}
For completeness, we include here a new proof of \eqref{OurPoissonSummation}, which is a variant of the Poisson summation formula.  This variant appears to be not very well known.  For a different proof, using locally compact abelian groups, see \cite{RichardStrungaru}. We also give proofs of some auxiliary lemmas that we require in the body of the paper, and which are also known. 
First, we recall  the following well-known Plancherel-Polya type inequality 
(see \cite{Schmeisse}{ page 216, Theorem 6}, for example).
\begin{lem} [Plancherel-Polya]
\label{Polya inequality}
Let $g:\R^d\rightarrow \C$ be a compactly supported function, with 
$\hat g \in L^1(\R^d)$.  Then there is a positive constant $M$ such that
\begin{equation}
\sum_{m\in\Lat}|\widehat{g}(m)|
\leq
M \int_{\mathbb{R}^{d}}|\widehat{g}(\xi)|.
\end{equation}
\end{lem}
\begin{proof}
Let $\psi$ be a  smooth and compactly supported function such that $\psi(x)=1$ for all $x$ in the support of $g$. So we have $g(x)=\psi(x)g(x)$ for all $x\in\R^d$. 
By our assumption that $\hat g \in L^1(\R^d)$ (and clearly $\hat \psi \in L^1(\R^d)$) we also have  $\widehat{g}(\xi)=(\widehat{\psi}*\widehat{g})(\xi)$.
We can now bound the series:
\begin{align}
\sum_{m\in\Lat}|\widehat{g}(m)|
&=\sum_{m\in\Lat}\left|\int_{\mathbb{R}^{d}} \widehat{\psi}(m-\xi) \widehat{g}(\xi)\,d\xi\right| \\
& \leq \sum_{m\in\Lat} \int_{\mathbb{R}^{d}}|\widehat{\psi}(m-\xi) \widehat{g}(\xi)|\,d \xi \\ \label{sum-integral interchange}
&=\int_{\mathbb{R}^{d}} \sum_{m\in\Lat}|\widehat{\psi}(m-\xi)||\widehat{g}(\xi)|\,d \xi \\
& \leq \sup _{\xi \in \mathbb{R}^{d}}\left(\sum_{m\in\Lat}|\widehat{\psi}(m-\xi)|\right) \int_{\mathbb{R}^{d}}|\widehat{g}(\xi)|\,d\xi \\
& = M \int_{\mathbb{R}^{d}}|\widehat{g}(\xi)|\,d\xi<\infty.
\end{align}
The integral-sum interchange in \eqref{sum-integral interchange} is justified by the monotone convergence Theorem. We observe that $\sum_{m\in\Lat}|\widehat{\psi}(m-\xi)|:= h(\xi)$ is a periodic function which is therefore completely determined on $\R^d/ \Lat$.  Due to the rapid decay of $\widehat{\psi}$, $h(\xi)$ is uniformly convergent on $\R^d/ \Lat$, and being a series of continuous functions, $h(\xi)$ is itself continuous. 
Since $\R^d/ \Lat$ is compact, $h$ attains its supremum there,  a finite constant $M>0$. 
\end{proof}
\begin{proof}
\hypertarget{Proof of Our Poisson Summation}
{(A Poisson Summation Formula, equation \eqref{OurPoissonSummation})}
The hypothesis that both $g$ and $\hat g \in L^1(\R^d)$ implies that we can use Fourier inversion (see \cite{Manfred}, Theorem 9.36). So $g(x) = 
\mathcal F({\hat g})(-x)$ is the image of an $L^1$ function under $\mathcal F$, and therefore uniformly continuous. Similarly $\hat g$ is uniformly continuous.  
Our goal is to prove that:
\begin{equation}\label{Poisson1}
    \sum_{n \in \Lat} g(n+x)=  \frac{1}{\det \Lat} \sum_{\xi \in \Lat^*} \hat{g}(\xi)e^{2\pi i \langle \xi, x\rangle}.
    \end{equation}
To this end, we first pick an infinitely smooth approximate identity  $\varphi_\varepsilon$
that is supported on the unit ball,
with $\varphi_\varepsilon \geq 0$ and $\int_{\R^d}\varphi_\varepsilon(x)\,dx=1$ for all $\varepsilon >0$ (i.e. $\varphi_\varepsilon$ is a bump function \cite{Shakarchi}, page 209).

Since $\varphi_{\varepsilon}$ vanishes outside a ball of radius $\varepsilon$ and $g$ is compactly supported, we have that $\varphi_{\varepsilon}*g$ is compactly supported on a  set that is independent of $\varepsilon$.  Therefore the sum 
\begin{equation*}
    \sum_{n\in\Lat}\varphi_\varepsilon*g (n+x):= G(x)
\end{equation*}
has a finite number of terms and is $\Lat$-periodic. Consequently, for
$m \in \Lat^*$ we have:

\begin{align*}
    \hat G(m) &:=\int_{{\R^d}/\Lat}\left(\sum_{n\in\Lat}\varphi_\varepsilon*g(n+x)\right)e^{-2\pi i \langle m,x\rangle}\,dx\\
    &=\sum_{n\in\Lat}\int_{{\R^d}/\Lat}\varphi_\varepsilon*g(n+x)e
    ^{-2\pi i \langle m,x\rangle}\,dx,
\end{align*}
so that we have:
\begin{align*}
    \hat G(m)
    &=\sum_{n\in\Lat}\int_{n+{\R^d}/\Lat}\varphi_\varepsilon*g(y)e^{-2\pi i \langle m,y-n\rangle}\,dy\\
    &=\sum_{n\in\Lat}\int_{n+{\R^d}/\Lat}\varphi_\varepsilon*g(y)e^{-2\pi i \langle m,y\rangle}\,dy\\
    &=\int_{\R^d}\varphi_\varepsilon*g(y)e^{-2\pi i \langle m,y\rangle}\,dy\\
    &=\widehat{\varphi_\varepsilon*g}(y)(m)
    =\widehat{\varphi}_\varepsilon(m)\widehat{g}(m)
    =\widehat{\varphi}(\varepsilon m)\widehat{g}(m).
\end{align*}
We notice that due to the fact that both
$\varphi_\varepsilon$ and $g$ are compactly supported, and that $\varphi_\varepsilon$ is infinitely smooth, 
$\varphi_\varepsilon*g$ belongs to the Schwartz class $S(\R^d)$. Thus, the basic Poisson summation formula for Schwartz functions holds: 
\begin{equation}
    \sum_{n\in\Lat}\varphi_\varepsilon*g(n+x)=\frac{1}{\det{\Lat}}\sum_{m\in\Lat^*}\widehat{\varphi}(\varepsilon m)\widehat{g}(m)e^{2\pi i \langle m, x\rangle},
\end{equation}
for every $x\in\R^d$. Because  $g$ is continuous and locally integrable on $\R^d$, the Lebesgue set of $g$ is $\R^d$ and we have
\begin{equation}
    \lim_{\varepsilon\rightarrow0^+}\sum_{n\in\Lat}\varphi_\varepsilon*g(n+x)=\sum_{n\in\Lat}\lim_{\varepsilon\rightarrow0^+}\varphi_\varepsilon*g(n+x)=\sum_{n\in\Lat}g(n+x).
\end{equation}
for every $x\in\R^d$.
Therefore 
\begin{align}\label{subtle limit}
    \sum_{n\in\Lat}g(n+x)
&=\lim_{\varepsilon\rightarrow0^+}\frac{1}{\det{\Lat}}\sum_{m\in\Lat^*}\widehat{\varphi}(\varepsilon m)\widehat{g}(m)e^{2\pi i \langle m, x\rangle}. 
\end{align}
Our next goal is to allow the interchange of the limit and the latter series in \eqref{subtle limit}, which is a subtle point. 
First, we have
\begin{equation}
    |\widehat{\varphi}(\varepsilon m)|=\left|\int_{\R^d}\varphi(x)e^{-2\pi i \langle x,\varepsilon m \rangle}\,dx\right|\leq\int_{\R^d}|\varphi(x)|\,dx=1.
\end{equation}
\noindent 
Summarizing, we have
$| \widehat{\varphi}(\varepsilon x)\widehat{g}(x)|
\leq |\widehat{g}(x)|$,
which is an absolutely summable dominating function because
\begin{equation}
\sum_{m\in\Lat^*}
 |\widehat{g}(m)|
 \leq 
M \int_{\mathbb{R}^{d}}|\widehat{g}(\xi)|
< \infty,
\end{equation}
for some positive constant $M$, by invoking 
Lemma \ref{Polya inequality} (a Plancherel-Polya inequality) with $\Lat$ replaced by
$\Lat^*$. 
By the Lebesgue dominated convergence theorem (applied to the counting measure on $\Lat^*$), we have
\begin{align*}
    \lim_{\varepsilon\rightarrow0^+}
    \sum_{m\in\Lat^*}
  \widehat{\varphi}(\varepsilon m)\widehat{g}(m)
     e^{2\pi i \langle m, x\rangle}
&=\sum_{m\in\Lat^*}\lim_{\varepsilon\rightarrow0^+}\widehat{\varphi}(\varepsilon m)\widehat{g}(m)e^{2\pi i \langle m, x\rangle} \\
&=\sum_{m\in\Lat^*}\widehat{g}(m)e^{2\pi i \langle m, x\rangle},
\end{align*}
where we've used the continuity of $\widehat\varphi$ in the last equality, and also 
$\widehat \varphi(0) = 1$.  Finally, putting together \eqref{subtle limit} with the latter computation, we have
\begin{equation}\label{General Poisson}
    \sum_{n\in\Lat}g(n+x)=\frac{1}{\det{\Lat}}\sum_{m\in\Lat^*}\widehat{g}(m)e^{2\pi i \langle m, x\rangle}
\end{equation}
for all $x\in\R^d$.
\end{proof}

We also include the following $3$ well-known lemmas, because often the precise relationship between convexity (or lack thereof) of $\QQ$ and the difference body $\QQ-\QQ$ can be subtle. 
\begin{lem}
\label{Q-Q = 2Q}
Let $\QQ\subset \R^d$ be a centrally symmetric body.
\begin{enumerate}[(a)]
    \item $\QQ - \QQ \supseteq 2\QQ$.
    \item If we also assume that $\QQ$ is convex, then $\QQ-\QQ=2\QQ$.
\end{enumerate}
\end{lem}
\begin{proof} To show (a), 
let $x\in 2\QQ$. Then $x=2y$ with $y\in\QQ$, implying that  $x=y+y=y-(-y)\in\QQ-\QQ$.  Here we used the central symmetry of $\QQ$
by invoking $-y\in\QQ$. 

To prove (b), it now suffices to show that $\QQ - \QQ \subseteq 2\QQ$.
So letting $x\in\QQ-\QQ$, $x=y-z$ with $y,z\in\QQ$. We can rewrite $x=y+(-z)=\left(\frac{1}{2}y+\frac{1}{2}(-z)\right)+\left(\frac{1}{2}y+\frac{1}{2}(-z)\right)$. Since $\QQ$ is centrally symmetric we have $-z\in\QQ$, and now the convexity of $\QQ$ implies that $\frac{1}{2}y+\frac{1}{2}(-z)\in\QQ$.  Therefore $x\in2\QQ$.     
\end{proof}
We observe that convexity is essential in part (b), by considering the following 
counterexample.  Let
 $C:= [-2, -1] \cup [1, 2]$, a nonconvex set in $\R$.  Here $C$ is centrally symmetric, yet
$C - C = [-3, 3] 
\not=[-4, -2] \cup [2, 4]=2C$.

\begin{lem}\label{basiclemma}
Let $A, B \subset \R^d$ be any two sets, and fix any 
$x\in\R^d$. Then:
\begin{equation}\label{aL2}
A \cap \big( B + x\big) \not= \varnothing
\iff   x \in A- B.
\end{equation}
Consequently, if $\QQ$ is a convex centrally symmetric body, then
\begin{equation}\label{bL2}
\QQ \cap \big( \QQ + x\big) \not= \varnothing
\iff   x \in 2\QQ.
\end{equation}
\end{lem}
\begin{proof}
Let $y\in A \cap (B+x)$.   Then
$y= z$ and $y = w + x$, where $z\in A$, and 
$w \in B$.
This gives us $z = w + x$, so that 
$x = y -  w \in A - B$, proving \eqref{aL2}. The converse follows exactly the same logical steps.    
To prove condition \eqref{bL2},  
we first set $A = B = \QQ$ in \eqref{aL2}. Because $\QQ$ is now also convex and centrally symmetric, 
Lemma \eqref{Q-Q = 2Q} implies that 
$ \QQ -  \QQ=2\QQ$, proving \eqref{bL2}. 
\end{proof}

The following is a standard exercise, but we include it here because it is used in Corollary \ref{application: arithmetic combintorics 1}.
\begin{lem}\label{lem:vanishing of the FT of the cube}  
Consider the unit cube 
$\square:= \left[ \tfrac{1}{2}, -\tfrac{1}{2} \right]^d$.
Then for $\xi \in \R^d$, we have
\[
\hat 1_{\square}(\xi)=0 \iff \text{ at least one of the coordinates of $\xi$ is a nonzero integer}.
\]
\end{lem}
\begin{proof}
We compute:
\begin{align}\label{FT of the cube}
   \hat 1_{\square}(\xi) &= 
\int_{\square} e^{-2\pi i \langle \xi, x \rangle} dx 
= \prod_{k=1}^d \int_{-\frac{1}{2}}^{\frac{1}{2}} e^{-2\pi i \langle \xi_k, x_k \rangle} dx_k.
\end{align}
The $k$'th term in the latter product \eqref{FT of the cube}
is just the sinc function 
${\rm sinc}(x_k)= \frac{\sin(\pi \xi_k)}{\pi \xi_k}$, 
when $\xi_k \not=0$, and ${\rm sinc}(0):=1$.
So the $k$'th term in the product defined by
\eqref{FT of the cube}
vanishes if and only if 
$\xi_k \in \Z\setminus \{0\}$, and we're done. 
\end{proof}
In other words, the zero set of the Fourier transform of the unit
cube consists precisely of the discrete collection of hyperplanes, defined by $\xi_k = n$, for any nonzero 
integer $n$, and any index $1\leq k \leq d$.


\end{document}